\documentclass [12pt] {amsart}

\usepackage [T1] {fontenc}
\usepackage {amssymb}
\usepackage {amsmath,bm}
\usepackage {amsthm}
\usepackage{tikz}
\usepackage{geometry}

\usepackage{url}

 \usepackage[titletoc,toc,title]{appendix}

\usetikzlibrary{arrows,decorations.pathmorphing,backgrounds,positioning,fit,petri}

\usetikzlibrary{matrix}

\usepackage{enumitem}
\usepackage{hyperref}

\makeatletter
\def\@tocline#1#2#3#4#5#6#7{\relax
  \ifnum #1>\c@tocdepth 
  \else
    \par \addpenalty\@secpenalty\addvspace{#2}%
    \begingroup \hyphenpenalty\@M
    \@ifempty{#4}{%
      \@tempdima\csname r@tocindent\number#1\endcsname\relax
    }{%
      \@tempdima#4\relax
    }%
    \parindent\z@ \leftskip#3\relax \advance\leftskip\@tempdima\relax
    \rightskip\@pnumwidth plus4em \parfillskip-\@pnumwidth
    #5\leavevmode\hskip-\@tempdima
      \ifcase #1
       \or\or \hskip 1em \or \hskip 2em \else \hskip 3em \fi%
      #6\nobreak\relax
    \dotfill\hbox to\@pnumwidth{\@tocpagenum{#7}}\par
    \nobreak
    \endgroup
  \fi}
\makeatother

\def\MLine#1{\par\hspace*{-\leftmargin}\parbox{\textwidth}{\[#1\]}}

\DeclareMathOperator {\td} {td}

\DeclareMathOperator {\im} {Im}
\DeclareMathOperator {\pr} {pr}

\DeclareMathOperator {\Der} {Der}
\DeclareMathOperator {\DF} {DF}

\DeclareMathOperator {\SL} {SL}

\DeclareMathOperator {\GL} {GL}

\DeclareMathOperator {\D} {D}

\DeclareMathOperator {\SAtyp} {SAtyp}

\DeclareMathOperator {\s} {s}
\DeclareMathOperator {\alg} {alg}

\DeclareMathOperator {\Jac} {Jac}
\DeclareMathOperator {\rk} {rk}

\DeclareMathOperator {\Zcl} {Zcl}

\DeclareMathOperator {\Atyp} {Atyp}
\DeclareMathOperator {\Span} {span}

\DeclareMathOperator {\C} {\mathbb{C}}

\DeclareMathOperator {\h} {\mathbb{H}}

\DeclareMathOperator {\Q} {\mathbb{Q}}

\DeclareMathOperator {\Ann} {Ann}
\DeclareMathOperator {\ppr} {Pr}

\DeclareMathOperator {\rar} {\rightarrow}
\DeclareMathOperator {\seq} {\subseteq}

\DeclareMathAlphabet\urwscr{U}{urwchancal}{m}{n}%
\DeclareMathAlphabet\rsfscr{U}{rsfso}{m}{n}
\DeclareMathAlphabet\euscr{U}{eus}{m}{n}
\DeclareFontEncoding{LS2}{}{}
\DeclareFontSubstitution{LS2}{stix}{m}{n}
\DeclareMathAlphabet\stixcal{LS2}{stixcal}{m} {n}

\newcommand{\dd}{\mathop{}\!\mathrm{d}}

\theoremstyle {definition}
\newtheorem {definition}{Definition} [section]
\newtheorem {remark} [definition] {Remark}
\newtheorem {remarks} [definition] {Remarks}
\newtheorem{example} [definition] {Example}
\newtheorem {claim} {Claim}
\newtheorem* {notation} {Notation}

\theoremstyle {plain}

\newtheorem {lemma} [definition] {Lemma}
\newtheorem {theorem} [definition] {Theorem}
\newtheorem {proposition} [definition] {Proposition}
\newtheorem {corollary} [definition] {Corollary}

\newtheorem {conjecture} [definition] {Conjecture}

\theoremstyle {remark}

\calclayout

\begin {document}

\title{Weak Modular Zilber--Pink with Derivatives}

\author{Vahagn Aslanyan}
\address{Department of Mathematical Sciences, Carnegie Mellon University, Pittsburgh, PA 15213, USA}

\address{Institute of Mathematics of the National Academy of Sciences, Yerevan 0019, Armenia}

\address{\texttt{Current address:} School of Mathematics, University of East Anglia, Norwich, NR4 7TJ, UK}
\email{V.Aslanyan@uea.ac.uk}

\thanks{This work was done while I was a postdoctoral associate at Carnegie Mellon University, and some revisions were made at the University of East Anglia. Partially supported by EPSRC grant EP/S017313/1.}

\date{}



\vspace*{-2cm}

\maketitle

\begin{abstract}
In unpublished notes Pila proposed a Modular Zilber--Pink with Derivatives (MZPD) conjecture, which is a Zilber--Pink type statement for the modular $j$-function and its derivatives. In this article we define D-special varieties, then state and prove two functional (differential) analogues of the MZPD conjecture for those varieties. In particular, we prove a weak version of MZPD. As a special case of our results, we obtain a functional Modular Andr\'e--Oort with Derivatives statement. The main tools used in the paper come from (model theoretic) differential algebra and complex analytic geometry, and the Ax--Schanuel theorem for the $j$-function and its derivatives (established by Pila and Tsimerman) plays a crucial role in our proofs. In the proof of the second Zilber--Pink type theorem we also use an Existential Closedness statement for the differential equation of the $j$-function.
\end{abstract}

\tableofcontents


\section{Introduction}

\subsection{The Conjecture on Intersections with Tori}

Schanuel conjectured (see \cite[p. 30]{Lang-tr}) that for any $\mathbb{Q}$-linearly independent complex numbers $z_1,\ldots,z_n$
$$\td_\mathbb{Q}\mathbb{Q}(z_1,\ldots,z_n, e^{z_1}, \ldots, e^{z_n}) \geq n$$ where $\td$ stands for transcendence degree. This is considered out of reach now. In particular, it implies algebraic independence of $e$ and $\pi$ which is a long-standing open problem. Nevertheless, Zilber gave an interesting model theoretic approach to Schanuel's conjecture. He constructed algebraically closed fields equipped with a unary function, called \emph{pseudo-exponentiation}, sharing some basic properties of the complex exponential function and, most importantly, satisfying (the analogue of) Schanuel's conjecture (see \cite{Zilb-pseudoexp}). He then showed that there is a unique (up to isomorphism) pseudo-exponential field of cardinality $2^{\aleph_0}$ and conjectured that it is isomorphic to $\mathbb{C}_{\exp}:= (\mathbb{C};+,\cdot, \exp, 0, 1)$.

In his work Zilber also considered a uniform version of Schanuel's conjecture, and observed that to deduce it from Schanuel's conjecture one needs a finiteness statement. He formulated a diophantine conjecture which serves that purpose -- the \emph{Conjecture on Intersections with Tori}, or CIT for short (see \cite{Zilb-exp-sum-published,Kirby-Zilber-exp}). We will formulate it shortly but we need to give some definitions first.

Below $(C;+,\cdot, 0,1)$ is an algebraically closed field of characteristic zero. The reader may assume $C$ is the field of complex numbers $\mathbb{C}$. All varieties are defined over $C$ unless explicitly stated otherwise and will be identified with the sets of their $C$-points. In particular, the affine space $\mathbb{A}^n$ will be denoted by $C^n$.  


Given two varieties $V, W$ inside a third variety $U$, the expected dimension of $V\cap W$ is $\dim V + \dim W - \dim U$. The intersection of two varieties in ``general position'' has the expected dimension. The following classical theorem from algebraic geometry states that in smooth varieties the dimension of an intersection is never less than the expected dimension.\footnote{This is not true without the assumption of smoothness. See \cite[p. 201, Footnote 20]{Loj-analytic-geom}.}

\begin{theorem}[Dimension of intersection]\label{dimension-of-intersection}
Let $U$ be a smooth irreducible algebraic variety and $V,W \subseteq U$ be irreducible subvarieties. Then any non-empty irreducible component $X$ of the intersection $V \cap W$ satisfies 
$$\dim X \geq \dim V+ \dim W - \dim U.$$
\end{theorem}

The generalisation of this theorem to complex analytic sets is also true and will be used in the proof of one of our main results. See \cite[Chapter III, 4.6]{Loj-analytic-geom} for a proof.

\begin{definition}
Let $V,W$ be subvarieties of a (not necessarily smooth) variety $U$. A non-empty irreducible component $X$ of $V \cap W$ is said to be \emph{typical} in $U$ if $\dim X = \dim V+ \dim W - \dim U$ and \emph{atypical} if $\dim X > \dim V+ \dim W - \dim U$.
\end{definition}

Thus, atypical components have atypically large (larger than expected) dimension. 

\begin{definition}
An \emph{algebraic torus} is an irreducible algebraic subgroup of $(C^{\times})^n$ for some positive integer $n$, where $C^{\times}$ is the multiplicative group of $C$.
\end{definition}

A variety defined by equations of the form $$y_1^{m_1}\cdots y_n^{m_n} =1,$$ where $m_i \in \mathbb{Z}$, is a subgroup of $(C^{\times})^n$ and can be decomposed into a disjoint union of an algebraic torus (the connected component of the identity element) and its torsion cosets. Note that an algebraic torus in $(\C^{\times})^n$ is the image of a $\Q$-linear subspace of $\C^n$ under the exponential map. 

\begin{definition}
Let $V\subseteq (C^{\times})^n$ be an algebraic variety. A subvariety $X\subseteq V$ is \emph{atypical} if it is an atypical (in $(C^{\times})^n$) component of an intersection $V\cap T$ where $T\subseteq (C^{\times})^n$ is a torsion coset of a torus.
\end{definition}


\begin{conjecture}[CIT, \cite{Zilb-exp-sum-published}]\label{cit-intro}
Let $V \subseteq (\mathbb{C}^{\times})^n$ be an algebraic variety.\footnote{This conjecture and the Modular Zilber--Pink conjecture stated below are usually formulated for the complex field but they make sense for an arbitrary algebraically closed field. In particular, the weak versions that we consider later are in that generality.} Then there is a finite collection $\Sigma$ of torsion cosets of proper subtori of $(\mathbb{C}^{\times})^n$ such that every atypical subvariety of $V$ is contained in some $T\in \Sigma$.
\end{conjecture}

Here torsion cosets of algebraic tori are the \emph{special} varieties. The form of the CIT conjecture is quite general in the sense that once there is a well-defined notion of special varieties having certain properties, one can formulate an analogous conjecture (see \cite[Conjecture 2.19]{Zilb-special-Schanuel}). 

Bombieri, Masser, and Zannier \cite{Bom-Mas-Zan} independently gave an equivalent formulation of CIT. Pink proposed a similar and more general conjecture for mixed Shimura varieties, again independently \cite{Pink,Pink-2}. The general conjecture is now known as the Zilber--Pink conjecture. It generalises the Mordell--Lang and Andr\'{e}--Oort conjectures, and also CIT.  

The Zilber--Pink conjecture (as well as CIT) is an active field of research in number theory and model theory. Many special cases are known, e.g. Mordell--Lang (Faltings, Raynaud, Vojta, Hindry, McQuillan), Andr\'{e}--Oort for arbitrary products of modular curves (Pila, \cite{Pila-Andre-Oort}) and for $\mathcal{A}_g$ -- the moduli space of principally polarized abelian varieties (Tsimerman, \cite{Tsimerman-AO-A_g}). See also \cite{Zannier-book-unlikely,Habegger-Pila-o-min-certain,Daw-Ren} for various other results around this conjecture.

Zilber showed in \cite{Zilb-exp-sum-published} that a functional analogue of Schanuel's conjecture established by Ax in \cite{Ax}, often called the Ax--Schanuel theorem, can be used to prove a weak version of the CIT conjecture \cite{Zilb-exp-sum-published}. 

\begin{theorem}[Weak CIT, \cite{Zilb-exp-sum-published,Kirby-semiab,Bom-Mas-Zan}]\label{weak-CIT}
For every subvariety $V\subseteq (C^{\times})^n$ there is a finite collection $\Sigma$ of proper subtori of $(C^{\times})^n$ such that every atypical component of an intersection of $V$ with a coset of a torus is contained in a coset of some torus $T \in \Sigma$.
\end{theorem}

We will not work with tori in this paper, nor will we need the above results; we presented them here to motivate the results of the paper and to give a brief account of the Zilber--Pink conjecture.


\subsection{Modular Zilber--Pink}

In the modular setting special varieties are defined in terms of modular polynomials (see Section \ref{jBackground}).

\begin{definition}
A $j$-\emph{special} subvariety of $\mathbb{C}^n$ (or, more generally, $C^n$, with coordinates $\bar{y}$) is an irreducible component of a variety defined by modular equations, i.e. equations of the form $\Phi_N(y_i, y_k) = 0$ for some $1\leq i, k \leq n$ where $\Phi_N(X,Y)$ is a modular polynomial.\footnote{More precisely, these are the special subvarieties of $Y(1)^n$ where $Y(1)$ is the modular curve $\SL_2(\mathbb{Z})\setminus \mathbb{H}$, which is identified with the affine line $\mathbb{C}$.} If no coordinate is constant on a $j$-special variety then it is said to be \emph{strongly} $j$-special.
\end{definition}

We call these varieties $j$-special since they are the images of \emph{special} subvarieties of Cartesian powers of the complex upper half-plane (these are defined by geodesic equations given by linear fractional transformations with rational coefficients) under the $j$-function (see Subsection \ref{subsection-MZPD}), and modular polynomials determine the ``functional equations'' of the $j$-function (similar to the equation $e^{x+y}=e^x \cdot e^y$ for the exponential function). See Section \ref{jBackground} for details. 

Given this notion of special varieties, atypical subvarieties are defined as above, that is, for a variety $V\subseteq \mathbb{C}^n$ and a $j$-special variety $S\subseteq \mathbb{C}^n$, a component $X$ of the intersection $V \cap S$ is a $j$-\emph{atypical} subvariety of $V$ if $\dim X > \dim V + \dim S - n$. Then the following is a modular analogue of CIT.

\begin{conjecture}[Modular Zilber--Pink Conjecture]
Let $V \subseteq \C^n$ be an algebraic variety. Then there is a finite collection $\Sigma$ of proper $j$-special subvarieties of $\C^n$ such that every $j$-atypical subvariety of $V$ is contained in some $T\in \Sigma$.
\end{conjecture}

This conjecture can be formulated in an equivalent form as follows (see \cite[Conjecture 7.1]{Pila-Tsim-Ax-j}).

\begin{conjecture}[MZP]\label{modular-ZP-intro}
Every algebraic variety in $\mathbb{C}^n$ contains only finitely many maximal $j$-atypical subvarieties.
\end{conjecture}

The following is yet another equivalent form of the Modular Zilber--Pink conjecture.

\begin{conjecture}[MZP]\label{conj-mzp-atyp_j}
Let $V\subseteq \C^n$ be an algebraic variety, and let $\Atyp_j(V)$ be the union of all $j$-atypical subvarieties of $V$. Then $\Atyp_j(V)$ is contained in a finite union of proper $j$-special subvarieties of $\C^n$.
\end{conjecture}

Note that the analogues of Conjectures \ref{modular-ZP-intro} and \ref{conj-mzp-atyp_j} for torsion cosets of algebraic tori are equivalent to CIT.

Zilber's argument for deducing weak CIT from Ax's theorem is quite general and goes through in various settings provided there is an appropriate analogue of Ax's theorem. In particular, the Ax--Schanuel theorem for the $j$-function established by Pila and Tsimerman in \cite{Pila-Tsim-Ax-j} (see Section \ref{jAS-ZP}) can be used to prove a weak version of the Modular Zilber--Pink conjecture. Below a \emph{strongly $j$-atypical} subvariety is a $j$-atypical subvariety with no constant coordinates.

\begin{theorem}[Weak Modular Zilber--Pink, {\cite[Theorem 7.1]{Pila-Tsim-Ax-j}}]\label{intro-thm-WMZP}
Every algebraic variety in $C^n$ contains only finitely many maximal strongly $j$-atypical subvarieties.
\end{theorem}

This theorem was proven by Pila and Tsimerman in \cite{Pila-Tsim-Ax-j} using o-minimality. We give a differential algebraic proof in Section \ref{weak-MZP} which is a direct analogue of the proof of weak CIT (we will follow Kirby's adaptation of Zilber's proof; see \cite[Theorem 4.6]{Kirby-semiab}). 

\subsection{Modular Zilber--Pink with Derivatives}\label{subsection-MZPD}

The $j$-function satisfies a third order differential equation and the aforementioned Ax--Schanuel theorem of Pila and Tsimerman  incorporates $j$ and its first two derivatives. There is also a Schanuel-type conjecture for $j, j', j''$, and the Ax--Schanuel theorem can be seen as a functional analogue of that conjecture. In unpublished notes Pila formulated a Modular Zilber--Pink with Derivatives conjecture -- henceforth referred to as MZPD -- and, assuming it, showed that the Modular Schanuel conjecture with Derivatives (MSCD) implies a uniform version of itself. This is an analogue of Zilber's result on the uniform Schanuel conjecture and CIT.\footnote{The MSCD conjecture was formulated by Pila in the same notes.} We now define special varieties for $j$ and its derivatives and present the MZPD conjecture following Pila's notes. The appropriate definitions can also be found in \cite{Spence}, where Spence discusses the Modular Andr\'e--Oort with Derivatives conjecture, also due to Pila.




While the definition of $j$-special varieties can be given without mentioning $j$ at all (and working with modular polynomials instead), it is more convenient to work with the $j$-function and its derivatives to define the special varieties in this setting. 
Let $\mathbb{H}$ be the complex upper half-plane and let $j: \mathbb{H}\rightarrow \mathbb{C}$ be the modular $j$-function; it is an analytic function on $\mathbb{H}$ (see Section \ref{jBackground}). Define a function $J:\mathbb{H}\rightarrow \mathbb{C}^3$ by $$J: z\mapsto (j(z), j'(z), j''(z)).$$ We extend $J$ to $\mathbb{H}^n$ by defining
$$J: \bar{z}\mapsto (j(\bar{z}),j'(\bar{z}),j''(\bar{z}))$$ where $j^{(k)}(\bar{z}) = (j^{(k)}(z_1), \ldots, j^{(k)}(z_n))$ for $k=0,1,2$. Note that we consider only the first two derivatives of $j$, for the higher derivatives are algebraic over those.

Let $\GL_2^+(\mathbb{Q})$ be the group of $2 \times 2$ rational matrices with positive determinant. This group acts on $\mathbb{H}$ by linear fractional transformations.

\begin{definition}\label{defin-intro-J-special}
\begin{itemize}[leftmargin = 0.5cm]
    \item[] 
    
    \item A subvariety $U\subseteq \mathbb{H}^n$ (i.e. an intersection of $\mathbb{H}^n$ with some algebraic variety) is called $\mathbb{H}$-\emph{special} if it is defined by some equations of the form $z_i = g_{i,k}z_k,~ i\neq k$, with $g_{i,k} \in \GL_2^+(\mathbb{Q})$, and some equations of the form $z_i = \tau_i$ where $\tau_i \in \mathbb{H}$ is a quadratic number.
    
    \item  For an $\h$-special variety $U$ we denote by $\langle \langle U \rangle \rangle$ the Zariski closure of $J(U)$ over $\mathbb{Q}^{\alg}$.
    
    \item A $J$-\textit{special} subvariety of $\mathbb{C}^{3n}$ is a set of the form $\langle \langle U \rangle \rangle$ where $U$ is a special subvariety of $\mathbb{H}^n$.\footnote{The notation $\langle \langle U \rangle \rangle$ is due to Pila, and the terms $\mathbb{H}$-special, $j$-special and $J$-special are due to Spence \cite{Spence}.}
\end{itemize}

\end{definition}

\begin{remark}\label{remark-image-analytic}
Note that for an $\mathbb{H}$-special $U\subseteq \mathbb{H}^n$ the set $j(U)\subseteq \mathbb{C}^n$ is defined by modular equations and is irreducible (since $U$ is irreducible), therefore it is $j$-special. Similarly, $J(U)$ is an irreducible locally analytic set and hence so is its Zariski closure. Thus, $J$-special varieties are irreducible.\footnote{Strictly speaking, $J(U)$ may not be complex analytic as it is the image of an analytic set under an analytic function, but it is locally analytic. It is irreducible in the sense that if $J(U)$ is contained in a union of analytic sets then it must be contained in one of them.}
\end{remark}


\begin{definition}
For a variety $V \subseteq \mathbb{C}^{3n}$ we let the $J$-\emph{atypical set} of $V$, denoted $\Atyp_J(V)$, be the union of all atypical components of intersections $V \cap T$ in $\mathbb{C}^{3n}$ where $T\subseteq \mathbb{C}^{3n}$ is a $J$-special variety.
\end{definition}

\begin{conjecture}[Pila, ``MZPD'']\label{MZPD-intro}
For every algebraic variety $V \subseteq \mathbb{C}^{3n}$ there is a finite collection $\Sigma$ of proper $\mathbb{H}$-special subvarieties of $\mathbb{H}^n$ such that $$\Atyp_J(V) \cap J(\mathbb{H}^n) \subseteq \bigcup_{\substack{U \in \Sigma\\ \bar{\gamma}\in \SL_2(\mathbb{Z})^n}} \langle \langle \bar{\gamma} U \rangle \rangle.$$ 
\end{conjecture}


\begin{remark}
Note that here we may need infinitely many $J$-special subvarieties to cover the atypical set of $V$ (see \cite{Spence}, discussion after Definition 1.5, and also Example \ref{example-AO-infinite} of this paper). Thus, according to this conjecture, the atypical set of $V$ is controlled by finitely many $\mathbb{H}$-special (equivalently, $j$-special) varieties but not necessarily finitely many $J$-special varieties.
\end{remark}

\subsection{Main results}

In this paper we explore functional analogues of Conjecture \ref{MZPD-intro}. Note that in modular Zilber--Pink without derivatives, and indeed in many other settings, a functional analogue of the main conjecture is also a weak version; in particular, special varieties are the same in both settings. However, the variants of the MZPD conjecture that we consider are not special cases, but rather some differential (functional) statements, and the special varieties that we work with are more general than $J$-special varieties. Nevertheless, our results do imply a weak version of the MZPD conjecture, which is formulated below. See Section \ref{section-analytic-MZPD} for details.

\begin{definition}
For a $J$-special variety $T\subseteq \mathbb{C}^{3n}$ and an algebraic variety $V \subseteq \mathbb{C}^{3n}$ an atypical component $X$ of an intersection $V \cap T$ in $\mathbb{C}^{3n}$ is a \emph{strongly $J$-atypical} subvariety of $V$ if for every irreducible analytic component $Y$ of $X \cap J(\mathbb{H}^n)$, no coordinate is constant on the projection of $Y$ on the first $n$ coordinates (corresponding to $j$-coordinates in $J$-special varieties). The \emph{strongly $J$-atypical set} of $V$, denoted $\SAtyp_J(V)$, is the union of all strongly $J$-atypical subvarieties of $V$.
\end{definition}

\begin{theorem}[Weak MZPD]\label{weak-MZPD1-intro}
For every algebraic variety $V \subseteq \mathbb{C}^{3n}$ there is a finite collection $\Sigma$ of proper $\mathbb{H}$-special subvarieties of $\mathbb{H}^n$ such that 
$$\SAtyp_J(V) \cap J(\mathbb{H}^n) \subseteq \bigcup_{\substack{U \in \Sigma\\ \bar{\gamma}\in \SL_2(\mathbb{Z})^n}} \langle \langle \bar{\gamma} U \rangle \rangle.$$ 
\end{theorem}

As pointed out above, the main results of this paper are some functional/differential analogues of the MZPD conjecture, where we deal with more general special varieties than $J$-special varieties. We call them D-\emph{special} varieties. 
We now define these varieties informally; a rigorous definition and an analysis of their structure is given in Section \ref{D-special}. 

Given a strongly $j$-special variety $T \subseteq C^n$ (with coordinates $\bar{y}$), i.e. a $j$-special variety with no constant coordinates, we consider a $C$-\emph{geodesic} subvariety $U$ of $C^n$ (with coordinates $\bar{x}$) defined as follows. For any two coordinates $y_i$ and $y_k$ which are related by a modular equation on $T$, i.e. $\Phi_N(y_i,y_k)=0$, we choose an arbitrary matrix $g_{i,k}\in \GL_2(C)$ and consider the equation $x_i = g_{i,k} x_k$ where $g_{i,k}$ is identified with the corresponding linear fractional transformation. 
These equations define $U$, which is called a $C$-geodesic variety associated with $T$. Note that there are infinitely many such varieties associated with the same $T$ since the $g_{i,k}$ are chosen arbitrarily.

Now let $(K;+, \cdot, D)$ be a differentially closed field with constant field $C$. Consider the set $W$ of all non-constant points $(\bar{z},\bar{j}, \bar{j}', \bar{j}'') \in U(K) \times T(K) \times K^{2n}$ which satisfy the differential equation of $j$, that is, $ Dj_i = j_i' D z_i,~ Dj'_i = j_i'' D z_i,~ Dj''_i = j_i''' D z_i$ where $j_i'''$ is determined from the equation $\Psi(j_i,j_i',j_i'',j_i''')=0$ where $\Psi$ is a rational function corresponding to the differential equation of $j$ as defined in Section \ref{jBackground}.\footnote{Thus, the function $\Psi(j(z),j'(z),j''(z),j'''(z))$ is identically zero where $j',j'',j'''$ are the complex derivatives of $j(z)$.} Here we think of each $z_i$ as some complex function and of $j_i, j_i', j_i'', j_i'''$ as $j(z_i), \frac{\dd j(z_i)}{\dd z_i}, \frac{\dd^2 j(z_i)}{\dd z_i^2}, \frac{\dd^3 j(z_i)}{\dd z_i^3}$ respectively, although complex solutions of the differential equation of the $j$-function are slightly more general (see Lemma \ref{solutions-of-eq-j}). Now consider the projection of $W$ onto the last $3n$ coordinates\footnote{The set $W$ can be thought of as a differential version of the set $\{ (\bar{z},J(\bar{z})): \bar{z}\in U \}$ from Definition \ref{defin-intro-J-special}. Then this projection corresponds to $J(U)$.}, i.e. onto $T(K)\times K^{2n}$, and take its Zariski closure over $C$. It is called a D-\emph{special} variety associated with $T$ and $U$ (it does not depend on the ambient differential field $K$). A D-special variety associated with $T$ is a D-special variety associated with $T$ and some $U$ (which must be associated with $T$). Note that D-special varieties are irreducible and, by definition, do not have constant coordinates.

This notion is the functional analogue of $J$-special varieties, and it is independent of the choice of the differentially closed field $K$. The equations defining D-special varieties are analysed in Section \ref{D-special}. One can use those equations to define D-special varieties without referring to a differential field. In Section \ref{section-analytic-MZPD} we generalise the notion of a $J$-special variety and define $J_{\bar{g}}$-special varieties where $\bar{g} \in \GL_2(\mathbb{C})^n$. They are defined like $J$-special varieties but with a function $J_{\bar{g}}$ instead of $J$, which is in some sense the composition of $J$ with $\bar{g}$ (where $\bar{g}$ is thought of as a vector function given by linear fractional transformations) and corresponds to a solution of the differential equation of the $j$-function. Then we show that $J_{\bar{g}}$-special varieties (for all $\bar{g}$) without constant coordinates are exactly D-special varieties over complex numbers (in particular, strongly $J$-special varieties are D-special). So, in the complex setting we have an analytic definition of D-special varieties.



We formulate a differential analogue of the MZPD conjecture below, but we fix some notation first. A D-\emph{atypical} subvariety of some variety $V\subseteq C^{3n}$ is an atypical (in $C^{3n}$) component of an intersection of $V$ with a D-special variety, and $\Atyp_{\D}(V)$ is the union of all D-atypical subvarieties of $V$. Further, in a differential field $(K; D)$ with constant field $C$, $E_J^{\times}(K)$ is the set of all non-constant tuples $(\bar{j}, \bar{j}', \bar{j}'')$
such that each $j_i$ satisfies the differential equation of the $j$-function for some $z_i$ and $j_i', j_i''$ are the derivatives of $j_i$ ``with respect to $z_i$''. More precisely, $(\bar{j}, \bar{j}', \bar{j}'')\in E_J^{\times}(K)$ if and only if for every $i$ we have
$$K \models \exists z_i, j_i''' ~ \left(  \Psi\left(j_i,j_i', j_i'', j_i'''\right)=0 \wedge  ( Dj_i=j_i'Dz_i \wedge Dj_i'=j_i''Dz_i\wedge Dj_i''=j_i'''Dz_i ) \right).$$

\begin{theorem}[Differential MZPD]\label{DMZPD-intro}
Let $(K; +,\cdot, D)$ be a differential field with an algebraically closed field of constants $C$. Given an algebraic variety $V\subseteq C^{3n}$, there is a finite collection $\Sigma$ of proper $j$-special subvarieties of $C^n$ such that $\Atyp_{\D}(V)(K) \cap E_J^{\times}(K)$ is contained in the union of all D-special varieties associated with $j$-special vareities from $\Sigma$.
\end{theorem}

The proof of this theorem is based on a uniform version of the Ax--Schanuel theorem for the $j$-function, Seidenberg's embedding theorem and some basic facts from complex analytic geometry such as Theorem \ref{dimension-of-intersection} for complex analytic sets. It is inspired by the proof of \cite[Theorem 7.1]{Pila-Tsim-Ax-j}, although our methods significantly differ from theirs and, in particular, we do not use o-minimality. In fact, we also give a second proof which is purely differential algebraic (modulo Ax--Schanuel).

Note that Pila and Scanlon have proven some differential algebraic Zilber--Pink theorems, without considering derivatives though.  Their results are similar to the above theorem in nature but the methods are quite different. They also work in a differential field and make implicit use of uniformity of differential Ax--Schanuel, although in a significantly different way. Of course, the most important difference is that we work with $j$ and its derivatives while they do not. We refer the reader to Scanlon's slides \cite{Pila-Scanlon-dif-ZP} for details.

Further, we also establish a ``more algebraic'' analogue of the MZPD conjecture.

\begin{theorem}[Functional MZPD]\label{FMZPD-intro}
Let $C$ be an algebraically closed field of characteristic zero. Given an algebraic variety $V\subseteq C^{3n}$, there is a finite collection $\Sigma$ of proper $j$-special subvarieties of $C^n$ such that every \emph{strongly} D-atypical subvariety of $V$ is contained in a D-special variety associated with some $T \in \Sigma$.
\end{theorem}

Note that here the definition of a \emph{strongly D-atypical} subvariety is more delicate than that of a strongly $j$-atypical subvariety. Roughly speaking, a strongly atypical subvariety is an atypical subvariety which is large enough to contain points of $E_J^{\times}(K)$ for a suitable differential field $K$. This corresponds to the intersection of $\Atyp_J(V)$ (or $\Atyp_{\D}(V)$) with the image of $J$ (respectively $E_J^{\times}$) in Conjecture  \ref{MZPD-intro} (Theorem \ref{DMZPD-intro}). For this reason, apart from the Ax--Schanuel theorem we also need an \emph{Existential Closedness} statement for the differential equation of the $j$-function (which was proposed as a conjecture in \cite{Aslanyan-adequate-predim} and was proven recently in \cite{Aslanyan-Eterovic-Kirby-Diff-EC-j}) to prove this theorem. Existential Closedness guarantees the existence of $E_J^{\times}$-points in strongly D-atypical subvarieties. 


As already mentioned above, we give a complex geometric proof for Theorem \ref{DMZPD-intro}, and a differential algebraic proof for Theorem \ref{FMZPD-intro}. However, we also show that the former can be deduced from the latter. Thus, we in fact give two proofs for Theorem \ref{DMZPD-intro}, one complex geometric and one differential algebraic.



Our results imply a functional version of the Modular Andr\'{e}--Oort with Derivatives (MAOD) conjecture (also proposed by Pila), which is an analogue of Andr\'e--Oort for the function $J$ and is a special case of MZPD (see Section \ref{section-FMAOD}).

\begin{theorem}[Functional MAOD]
Let $C$ be an algebraically closed field of characteristic zero. Given an algebraic variety $V\subsetneq C^{3n}$, there is a finite collection $\Sigma$ of proper $j$-special subvarieties of $C^n$ such that every D-special subvariety of $V$ is contained in a D-special variety associated with some $T \in \Sigma$.
\end{theorem}

Note that Spence also proved a weak version of the MAOD conjecture using o-minimality techniques in \cite{Spence}. It seems that his and our results are two different statements related to the MAOD conjecture and none of them follows from the other. Spence's version is actually a weakening of MAOD where he deals with $J$-special subvarieties which may not be strongly special, and is a number theoretic result, while our version is indeed a functional analogue of MAOD (and we do not deal with varieties with constant coordinates). Spence proved the full MAOD conjecture assuming a transcendence conjecture for the values of $j'$.

Note also that we consider uniform versions of all the aforementioned results and prove them for parametric families of varieties.


Although all results of the paper are related to each other, the reader may skip some sections if they are interested only in specific theorems. For instance, for a proof of Theorem \ref{weak-MZPD1-intro} one may read Sections \ref{jBackground}, \ref{jAS-ZP} and \ref{section-analytic-MZPD} only.

\addtocontents{toc}{\protect\setcounter{tocdepth}{1}}
\subsection*{Notation and conventions}

We fix some notation and conventions here that will be used throughout the paper.

\begin{itemize}[leftmargin=0.5cm]

    \item For a set $A$ and a tuple $\bar{a}\in A^m$ we will  write $\bar{a} \subseteq A$ when the length of $\bar{a}$ is not essential.
    
    \item All fields considered in this paper are of characteristic $0$.
    
    \item 
    Irreducible varieties are assumed to be absolutely irreducible. 
    
    \item When we work over some algebraically closed field $C$, we will identify algebraic varieties defined over $C$ with the sets of their $C$-rational points and write $V\subseteq C^n$. For a larger field $F\supseteq C$ the set of $F$-points of $V$ is denoted by $V(F)$. 
    
    
     \item By ``generic'' we always mean generic in the sense of fields, i.e. Zariski generic, unless explicitly stated otherwise.

    
    \item If $W\subseteq C^n \times C^m$ is an algebraic variety then for a tuple $\bar{c} \in C^m$ we denote by $W_{\bar{c}}$ the projection on $C^n$ of the fibre of the projection map $\pr: W\rightarrow C^m$ above $\bar{c}$, that is, $W_{\bar{c}}(C) = \{ \bar{w}\in C^n: (\bar{w},\bar{c})\in W(C)\}$. Then the collection $(W_{\bar{c}})_{\bar{c} \subseteq C}$ is a parametric family of varieties in $C^n$. This notation is not actually precise since we should let $\bar{c}$ vary over the $C$-points of the projection of $W$, but we will write $\bar{c} \subseteq C$ for simplicity.
    
    \item If $K \subseteq F$ are fields, the transcendence degree of $F$ over $K$ will be denoted by $\td_KF$ or $\td(F/K)$. The algebraic closure of a field $F$ is denoted by $F^{\alg}$.
    
    \item For fields $K \subseteq F$ and a set $X \subseteq F^n$ the Zariski closure of $X$ over $K$ will be denoted $\Zcl(X/K)$.


    \item When we work in the affine space $F^{2n}$ (for some field $F$), we will denote its coordinates by $(x_1,\ldots,x_n, y_1,\ldots,y_n)$, or concisely $(\bar{x},\bar{y})$. The coordinates of $F^{4n}$ (or $F^{3n}$) will be denoted by $(\bar{x},\bar{y},\bar{y}',\bar{y}'')$ (respectively $(\bar{y},\bar{y}',\bar{y}'')$). We will use $\bar{x}$ or $\bar{y}$ (depending on the context) for the coordinates of $F^n$. 
    
    Note that $'$ is just a symbol here and not a derivation, although when we work in a differential field, $y_i'$ and $y_i''$ will correspond to derivations of $y_i$ with respect to $x_i$ in some precise sense which will be made clear later. We never denote a derivation of an abstract differential field by $'$. However, when we work with actual complex functions, we will use $'$ to denote their derivatives with respect to their arguments. We often state this explicitly to avoid any possible confusion.

    
    \item In a differential field $(K;+,\cdot, D)$ for a non-constant element $x\in K$ we define a derivation $\partial_x:K\rightarrow K$ by $\partial_x: y \mapsto \frac{Dy}{Dx}$ and often call it \emph{differentiation with respect to} $x$.
    
    \item $\Psi(Y_0,Y_1,Y_2,Y_3)$ is a rational function (over $\mathbb{Q}$) corresponding to the differential equation of the $j$-function, i.e. $\Psi(j, j', j'', j''') = 0$ where $j^{(k)} = j^{(k)}(z)$ is the $k$-th derivative of the $j$-function (see Section \ref{jBackground}).
    
\end{itemize}
\addtocontents{toc}{\protect\setcounter{tocdepth}{2}}



\section{Differential algebraic preliminaries}\label{section-dif-alg-prelim}
\setcounter{equation}{0}

We assume familiarity with the basics of differential algebra and model theory of differential fields. The reader is referred to \cite{Mar-dif} for an introduction to the topic. Nevertheless, in this section we introduce some preliminary concepts and results that will be used in the proofs of the main theorems of the paper. 

\subsection{Differential forms}\label{diff-forms}

Let $C\subseteq K$ be fields of characteristic zero. The vector space of \textit{abstract differential forms} (or \emph{K\"{a}hler differentials}) on $K$ over $C$, denoted $\Omega(K/C)$, is the quotient of the vector space generated by the set of symbols $\{ dx: x \in K \}$ by the relations
$$d(x+y)=dx+dy,~ d(xy)=xdy+ydx,~ dc=0,~ c\in C.$$ When no confusion can arise, we may drop $K/C$ from the notation of $\Omega$.

The map $d:K \rightarrow \Omega(K/C)$ is the \emph{universal derivation} on $K$. It satisfies the following universal property: for every $K$-vector space $V$ and every derivation $\delta: K \rightarrow V$ with $\delta|_C = 0$ there is a unique linear map $\xi: \Omega \rightarrow V$ such that $\xi \circ d = \delta$.

It is easy to verify that arbitrary elements $x_1, \ldots, x_n \in K$ are algebraically dependent over $C$ if and only if $dx_1, \ldots, dx_n$ are linearly dependent over $K$. Indeed, differentiating a relation $p(x_1,\ldots,x_n)=0$ we get a linear relation for $dx_1,\ldots,dx_n$. The converse follows from the universal property of $d$.

In particular, $\td(K/C) = \dim \Omega(K/C)$ and $\ker(d)$ is equal to the relative algebraic closure of $C$ in $K$. From now on we assume $C$ is relatively algebraically closed in $K$, therefore $\ker(d)=C$.

The vector space of derivations on $K$ that vanish on $C$ is denoted by $\Der(K/C)$. A differential form $\omega \in \Omega$ can be thought of as a linear functional $\omega: \Der(K/C) \rightarrow K$. For $x\in K$ we define $dx(D) = Dx$ for every $D \in \Der(K/C)$ and extend it to $\Omega$ by linearity. Thus, differential forms on $K$ over $C$ can be defined as linear forms on $\Der(K/C)$. This establishes an embedding of $\Omega$ into $(\Der(K/C))^*$, the dual space of $\Der(K/C)$. When $\td(K/C)$ is finite, we see that $\dim (\Der(K/C))^* = \dim \Der(K/C) = \td(K/C) = \dim \Omega(K/C)$ hence the above embedding is an isomorphism. So, $\Omega(K/C)$ can be identified with the space $(\Der(K/C))^*$. Furthermore, this shows that $\Der(K/C)$ can be identified with the dual $(\Omega(K/C))^*$. This identification can be described explicitly: a derivation $D\in \Der(K/C)$ is identified with its dual $$D^*: \Omega(K/C)\rightarrow K,~ D^*: dx \mapsto Dx.$$

Assume for some elements $x_1,\ldots,x_n\in K$ we have $C(\bar{x})\subseteq K\subseteq C(\bar{x})^{\alg}$. If $\td(K/C) = m$ then $\dim \Der(K/C) = m$. Let $D_1,\ldots,D_l\in \Der(K/C)$ be $K$-linearly independent derivations (where $l\leq m$). Consider the Jacobian matrix $\Jac(\bar{x}) := (D_ix_k)_{k,i}$. 

\begin{claim}
$\rk \Jac(\bar{x}) = l$.
\end{claim}
\begin{proof}
Assume $\rk \Jac(\bar{x}) < l$. Then there are $a_1,\ldots,a_l\in K$, not simultaneously zero, such that
$$\sum_{i=1}^l a_i \cdot D_i x_k =0, \mbox{ for all } k=1,\ldots,n.$$
Consider the derivation $D:=\sum_i a_i D_i$. Then for every $k$ we have $$D^*(dx_k)=0.$$ Since $C(\bar{x})\subseteq K\subseteq C(\bar{x})^{\alg}$ and $D^*$ is linear, we conclude that $D^*(\omega) = 0$ for any $\omega \in \Omega(K/C)$. Thus, $D^* = 0$, hence $D=0$ which means $D_1,\ldots,D_l$ are linearly dependent, a contradiction.
\end{proof}

\subsection{Non-commuting derivations}\label{non-commuting-der}

The Ax--Schanuel theorem for the $j$-function that we consider in Section \ref{jAS-ZP} holds for differential fields with commuting derivations. However, we need a slightly general version of the theorem where the derivations satisfy a weaker condition than commutativity. In this section we introduce the necessary tools for generalising some statements from commuting to non-commuting derivations.

Let $C\subseteq K$ be fields. We define the \emph{Lie bracket} on $\Der(K/C)$ by $$[D_1,D_2] = D_1 \circ D_2 - D_2 \circ D_1, \mbox{ i.e. } [D_1,D_2]x = D_1 (D_2x) - D_2(D_1x).$$ It is easy to verify that $[D_1,D_2] \in \Der(K/C)$ and $\Der(K/C)$ is a Lie algebra over $C$. Two derivations commute if and only if their Lie bracket is zero.

The following is an analogue of the Frobenius theorem of differential geometry. See \cite[Chapter 0, \S 5, Proposition 6]{Kolchin-diff-alg-gp} or \cite[Lemma 2.2]{Singer-noncommuting}  for a proof.

\begin{lemma}\label{Lie-Lemma}
Let $D_1,\ldots, D_m \in \Der(K)$ be linearly independent (over $K$) derivations. Assume that for each  $i, j$
\begin{equation}\label{Lie}
    [D_i, D_j] \in \Span_K \{ D_1, \ldots, D_m \}.
\end{equation}
Then there exist $K$-linearly independent commuting derivations $$\delta_1, \ldots, \delta_m \in \Span_K \{ D_1, \ldots, D_m \}.$$
\end{lemma}

In other words, any finite dimensional space of derivations which is closed under the Lie bracket has a commuting basis.

\begin{definition}
Let $(K;+,\cdot, D_1, \ldots, D_m)$ be a differential field with $m$ derivations. We say $K$ is a \emph{Lie differential field} iff the condition \eqref{Lie} is satisfied.
\end{definition}

Note that in this definition we do not assume that the derivations are linearly independent. However, if \eqref{Lie} holds for a $K$-linear basis of $\Span_K \{ D_1, \ldots, D_m \}$ then the latter is closed under the Lie bracket.

\section{Background on the $j$-function}\label{jBackground}
\setcounter{equation}{0}

We do not need to know much about the $j$-function itself, nor need we know its precise definition. Being familiar with some basic properties of $j$ will be enough for this paper. We summarise those properties below referring the reader to \cite{Lang-elliptic,Serre,Masser-Heights,Silverman} for details.

Let $\GL_2(\mathbb{C})$ be the group of $2 \times 2$ complex matrices with non-zero determinant. This group acts on on the Riemann sphere by linear fractional transformations. Namely, for a matrix $g=\begin{pmatrix}
a & b \\
c & d
\end{pmatrix} \in \GL_2(\mathbb{C})$ we define 
$$gz=\frac{az+b}{cz+d}.$$
This action is obviously the same as the action of the subgroup $\SL_2(\mathbb{C})$ consisting of matrices with determinant $1$ (to be more precise, the action of $\GL_2(\mathbb{C})$ factors through $\SL_2(\mathbb{C})$).

The function $j$ is a modular function of weight $0$ for the \emph{modular group} $\SL_2(\mathbb{Z})$, which is defined and analytic on the upper half-plane $\mathbb{H}:=\{ z \in \mathbb{C}: \im(z)>0 \}$. It is $\SL_2(\mathbb{Z})$-invariant. Moreover, by means of $j$ the quotient $Y(1):=\SL_2(\mathbb{Z}) \setminus \mathbb{H}$ is identified with $\mathbb{C}$ (thus, $j$ is a bijection from the fundamental domain of $\SL_2(\mathbb{Z})$ to $\mathbb{C}$). 

\subsection{Modular polynomials}
Let $\GL_2^+(\mathbb{R})$ be the subgroup of $\GL_2(\mathbb{R})$ consisting of matrices with positive determinant.\footnote{This group acts on the upper half-plane and in fact it is the largest subgroup of $\GL_2(\mathbb{C})$ with that property.} Let $\GL_2^+(\mathbb{Q})$ be its subgroup of matrices with rational entries. For $g \in \GL_2^+(\mathbb{Q})$ we let $N(g)$ be the determinant of $g$ scaled so that it has relatively prime integral entries. For each positive integer $N$ there is an irreducible polynomial $\Phi_N(X,Y)\in \mathbb{Z}[X,Y]$ such that whenever $g \in \GL_2^+(\mathbb{Q})$ with $N=N(g)$, the function $\Phi_N(j(z),j(gz))$ is identically zero. Conversely, if $\Phi_N(j(x),j(y))=0$ for some $x,y \in \mathbb{H}$ then $y=gx$ for some $g \in \GL_2^+(\mathbb{Q})$ with $N=N(g)$. The polynomials $\Phi_N$ are called \emph{modular polynomials}. It is well known that $\Phi_1(X,Y)=X-Y$ and all the other modular polynomials are symmetric. The properties of modular polynomials imply, in particular, that if $\tau\in \mathbb{H}$ is a quadratic number then $j(\tau)$ is algebraic. These numbers are known as \emph{special values} of $j$ or as \emph{singular moduli}. For $w=j(z)$ the image of the $\GL_2^+(\mathbb{Q})$-orbit of $z$ under $j$ is called the \emph{Hecke orbit} of $w$. It obviously consists of the union of solutions of the equations $\Phi_N(X,w)=0,~ N\geq 1$. Two elements $w_1,w_2 \in \mathbb{C}$ are called \emph{modularly independent} if they have different Hecke orbits, i.e. do not satisfy any modular relation $\Phi_N(w_1,w_2)=0$. This definition makes sense for arbitrary fields (of characteristic zero) as the modular polynomials have integer coefficients.

\subsection{Differential equation}
The $j$-function satisfies an order $3$ algebraic differential equation over $\mathbb{Q}$, and none of lower order (i.e. its differential order over $\mathbb{C}$ is $3$). Namely, $\Psi(j,j',j'',j''')=0$ where
$$\Psi(Y_0,Y_1,Y_2,Y_3)=\frac{Y_3}{Y_1}-\frac{3}{2}\left( \frac{Y_2}{Y_1} \right)^2 + \frac{Y_0^2-1968Y_0+2654208}{2Y_0^2(Y_0-1728)^2}\cdot Y_1^2.$$

Notice that
$$\Psi(Y,Y',Y'',Y''')=S(Y)+R(Y)(Y')^2,$$
where $S$ denotes the \emph{Schwarzian derivative} defined by $$S(Y) = \frac{Y'''}{Y'} - \frac{3}{2} \left( \frac{Y''}{Y'} \right) ^2,$$ and  $$R(Y)=\frac{Y^2-1968Y+2654208}{2Y^2(Y-1728)^2}.$$ \textit{Throughout the paper $\Psi$ will always denote the above rational function.} Observe also that $\Psi$ is linear with respect to $Y'''$ so the differential equation of the $j$-function can be written as $Y'''=\eta(Y,Y',Y'')$ for some rational function $\eta(Y_0,Y_1,Y_2)$.\footnote{Here we think of $Y$ as a variable ranging over the set of complex functions defined on some one-dimensional domain, and so $Y'$ denotes the derivative of such a function with respect to its argument. Later, when we work in an abstract differential field, we will interpret this differential equation by replacing $'$ with an appropriate derivation (which may not be the derivation of the field). Recall that we do not use $'$ to denote a derivation of an abstract differential field, and when $j$ is an element in such a field, $j', j'', j'''$ will be just some elements that will normally correspond to the abstract derivatives of the element $j$ with respect to some element $z$ (which may be thought of as the argument of $j$).}

The following result is well known (see, for example, \cite[Lemma 4.2]{Freitag-Scanlon} or \cite[Lemma 4.1]{Aslanyan-adequate-predim} for a proof).
\begin{lemma}\label{solutions-of-eq-j}
All functions $j(gz)$ with $g \in \SL_2(\mathbb{R})$ satisfy the differential equation $\Psi(y, y', y'', y''')=0$ and all solutions of  that equation defined on $\mathbb{H}$ are of that form. If we allow functions not necessarily defined on $\mathbb{H}$, then all solutions will be of the form $j(gz)$ where $g \in \SL_2(\mathbb{C})$.
\end{lemma}

\subsection{Functional equations of the differential equation of the $j$-function}

Let us introduce some notation first. Below $y',y'',y'''$ are variables and do not denote the derivatives of $y$.

\begin{notation} Let $(K;+,\cdot, D_1,\ldots,D_m)$ be a differential field with $m$ derivations and let $C:= \bigcap_i \ker D_i$ be the field of constants.
\begin{itemize}[leftmargin=0.5cm]
    \item Let $E_{(z,J)}(x,y,y',y'')$ denote the formula

        \MLine{\exists y''' ~ \left(  \Psi\left(y,y', y'', y'''\right)=0 \wedge  \bigwedge_{k=1}^m ( D_ky=y'D_kx \wedge D_ky'=y''D_kx\wedge D_ky''=y'''D_kx ) \right).}

By abuse of notation for any $n$ we will also let $E_{(z,J)}(K)$ denote the set of all tuples $(\bar{z},\bar{j},\bar{j}',\bar{j}'')\in K^{4n}$ with $(z_i,j_i,j_i',j_i'')\in E_{(z,J)}(K)$. The set $E^{\times}_{(z,J)}(K)$ consists of all $E_{(z,J)}(K)$-points that do not have any constant coordinates.

\item  $E_{(z,j)}(x,y)$ is the projection $\exists y', y'' E_{(z,J)}(x,y,y',y'')$. As above $E_{(z,j)}(K)$ also denotes the set of all tuples $(\bar{z},\bar{j})\in K^{2n}$ for which $(z_i,j_i)\in E_{(z,j)}(K)$, and $E^{\times}_{(z,j)}(K)$ consists of all $E_{(z,j)}(K)$-points that do not have any constant coordinates.

\item $E_J(y,y',y'')$ is the projection of $E_{(z,J)}$ onto the last three coordinates, i.e. $E_J$ denotes the formula $\exists x E_{(z,J)}(x,y,y',y'').$ Equivalently, $E_J$ is given by
\MLine{\exists y'''\left( \Psi(y,y',y'',y''')=0 \wedge \bigwedge_{k=1}^m \frac{D_ky}{y'} = \frac{D_ky'}{y''} = \frac{D_ky''}{y'''}\right).}
As above, $E_J(K)$ also denotes the set $\{ \bar{J} = (\bar{j},\bar{j}',\bar{j}''): (j_i, j_i', j_i'')\in E_J(K) \mbox{ for all } i\}$, and $E^{\times}_J(K)$ is the set of all points in $E_J(K)$ with no constant coordinates.
\end{itemize}
\end{notation}

The next proposition describes the ``functional equations'' of the relations $E_{(z,j)}$ and $E_{(z,J)}$.

\begin{proposition}[{\cite[Lemmas 4.10, 4.11 and 4.41]{Aslanyan-adequate-predim}}]\label{funct-eq}
Let $(F; +, \cdot, D)$ be a differential field with field of constants $C$. 
\begin{enumerate}
    \item[\textup{(i)}] If $(z_i,j_i) \in E^{\times}_{(z,j)}(F),~ i=1,2,$ and $\Phi_N(j_1, j_2)=0$ for some modular polynomial $\Phi_N$ then $z_2 = gz_1$ for some $g \in \SL_2(C)$.
    \item[\textup{(ii)}] If $(z_1,j_1)\in E^{\times}_{(z,j)}(F)$ and $(z_2,j_2)\in F^2$ such that $\Phi_N(j_1,j_2)=0$ for some $\Phi_N$ and $z_2=gz_1$ for some $g\in \SL_2(C)$ then $(z_2,j_2)\in E^{\times}_{(z,j)}(F)$.
    
    \item[\textup{(iii)}] If $\left(z,j,j',j''\right) \in E^{\times}_{(z,J)}(F)$ then for any  $g=\begin{pmatrix} a & b \\ c & d \end{pmatrix} \in \SL_2(C)$ 
$$\left(gz,j,j'\cdot (cz+d)^2,j''\cdot (cz+d)^2-2c\cdot j'\cdot (cz+d)^3\right) \in E^{\times}_{(z,J)}(F).$$
Conversely, if for some $j$ we have $\left(z_1,j,j',j''\right), \left(z_2,j, w',w''\right) \in E^{\times}_{(z,J)}(F)$ then $z_2=gz_1$ for some $g \in \SL_2(C)$.

\item[\textup{(iv)}] If $\left(z,j_1,j_1',j_1''\right) \in E^{\times}_{(z,J)}(F)$ and $\Phi(j_1,j_2)=0$ for some modular polynomial $\Phi(X,Y)$ then $\left(z,j_2,j_2',j_2''\right) \in E^{\times}_{(z,J)}(F)$ where $j_2',~ j_2''$ are determined from the following system of equations:
\begin{align*}
& \frac{\partial \Phi}{\partial X}(j_1,j_2)\cdot j_1' + \frac{\partial \Phi}{\partial Y}(j_1,j_2)\cdot j_2' =0,\\
& \begin{aligned}
\frac{\partial^2 \Phi}{\partial X^2}(j_1,j_2)\cdot \left(j_1'\right)^2+ \frac{\partial^2 \Phi}{\partial Y^2}(j_1,j_2)\cdot \left(j_2'\right)^2 + 2\cdot \frac{\partial^2 \Phi}{\partial X \partial Y}(j_1,j_2)\cdot j_1'\cdot j_2'+ & \\ 
+ \frac{\partial \Phi}{\partial X}(j_1,j_2)\cdot j_1'' + \frac{\partial \Phi}{\partial Y}(j_1,j_2) & \cdot j_2''=  0.
\end{aligned}
\end{align*}
\end{enumerate}
\end{proposition}

Note that (iv) above can be generalised to different $z_1$ and $z_2$ linked by an $\SL_2(C)$-relation. However, the general case follows from the above properties.

\begin{remark}
The converse of (i) is not true: if $z_2=gz_1$ for some $g\in \SL_2(C)$ then this does not impose a relation on $j_1,j_2$, they may be algebraically independent (see also Remark \ref{remark-alg-ind-sol}).
\end{remark}

\section{Ax--Schanuel for the $j$-function}\label{jAS-ZP}
\setcounter{equation}{0}

The following theorem was proved by Pila and Tsimerman in \cite{Pila-Tsim-Ax-j}.

\begin{theorem}[Ax--Schanuel with Derivatives for $j$]\label{j-chapter-Ax-for-j}
Let $(K; +, \cdot, D_1,\ldots,D_m)$ be a differential field with commuting derivations and with field of constants $C$. Assume $(z_i, j_i, j_i', j_i'') \in E^{\times}_{(z,J)}(K),~ i=1,\ldots,n$.
If the $j_i$'s are pairwise modularly independent then 
\begin{equation}\label{j-chapter-Ax-ineq}
\td_CC\left(\bar{z},\bar{j},\bar{j}',\bar{j}''\right) \geq 3n+\rk \Jac(\bar{z}).
\end{equation}
\end{theorem}

\begin{lemma}\label{Ax-for-Lie}
The Ax--Schanuel theorem for the $j$-function holds in all Lie differential fields.
\end{lemma}
\begin{proof}
Let $(K;+,\cdot, D_1,\ldots,D_m)$ be a Lie differential field and let $z_i, j_i, j_i', j_i'', j_i'''$ be as in Theorem \ref{j-chapter-Ax-for-j}. If we replace $D_1,\ldots,D_m$ by a basis of $\Span_K\{ D_1, \ldots, D_m \}$ then it will not affect the Jacobian of any tuple. Hence we may assume that $D_1,\ldots,D_m$ are linearly independent over $K$. By Lemma \ref{Lie-Lemma} there are $K$-linearly independent and commuting derivations $$\delta_1, \ldots, \delta_m \in \Span_K \{ D_1, \ldots, D_m \}.$$ But then it is clear that $$\bigcap_i \ker D_i = \bigcap_i \ker \delta_i$$ and 
\begin{equation*}
\Psi\left(j_i,j_i', j_i'', j_i'''\right)=0 \wedge \delta_kj_i=j_i'\delta_k z_i \wedge \delta_kj_i'=j_i''\delta_kz_i\wedge \delta_kj_i''=j_i'''\delta_kz_i
\end{equation*}
for all $k=1,\ldots,m$. Hence the inequality \eqref{j-chapter-Ax-ineq} holds and $\rk (D_i x_k)_{i,k} = \rk (\delta_i x_k)_{i,k}$.
\end{proof}

Note however that in applications we find it more convenient to use Lemma \ref{Lie-Lemma} to choose commuting derivations of the field under consideration and apply Theorem \ref{j-chapter-Ax-for-j} instead of applying Lemma \ref{Ax-for-Lie} directly. 

The following statement is a direct consequence of Theorem \ref{j-chapter-Ax-for-j} (see also \cite[\S 2.6]{Pila-Tsim-Ax-j}).

\begin{theorem}[Ax--Schanuel without derivatives]\label{j-Ax-without-der}
Let $(K; D_1,\ldots,D_m)$ be a differential field with commuting derivations and with field of constants $C$. Assume $(z_i, j_i) \in E^{\times}_{(z,j)}(K),~ i=1,\ldots,n$.
If the $j_i$'s are pairwise modularly independent then 
\begin{equation*}
\td_CC\left(\bar{z},\bar{j}\right) \geq n+\rk \Jac(\bar{z}).
\end{equation*}
\end{theorem}

This theorem implies in particular that the only algebraic relations between the functions $j(z)$ and $j(gz)$ for $g \in \GL_2^+(\mathbb{R})$ are the modular relations (corresponding to $g \in \GL^+_2(\mathbb{Q})$). 

\begin{remark}\label{remark-alg-ind-sol}
As pointed out above, if $E^{\times}_{(z,j)}(z_i,j_i),~ i=1,2$, and $j_1, j_2$ are modularly dependent then $z_1$ and $z_2$ are $\SL_2(C)$-related, but if $z_1$ and $z_2$ are $\SL_2(C)$-related then $j_1$ and $j_2$ may still be algebraically independent over $C$. Nevertheless, in that case we know by Ax--Schanuel that $j_1$ and $j_2$ must be either algebraically independent (over $C(z_1)$) or related by a modular relation (see also \cite{Aslanyan-SM-reducts}).
\end{remark}

Let us establish uniform versions of the above theorems. 

\begin{theorem}[Uniform Modular Ax--Schanuel with Derivatives]\label{uniform-AS-j}
Let $S\seq (\Q^{\alg})^n$ be a $j$-special variety and let $W\seq (\Q^{\alg})^n \times S(\Q^{\alg})\times (\Q^{\alg})^{2n} \times (\Q^{\alg})^l$ be an algebraic variety defined over $\Q^{\alg}$. Then there is a finite  set $\Sigma = \Sigma(W)$ (depending on $W$) of proper $j$-special subvarieties of $S$ with the following property.

For every differential field $(K; D_1,\ldots,D_m)$ with $m$ commuting derivations and with field of constants $C$, for every  $\bar{c} \in C^l$ and for every\footnote{Recall that for $\bar{c}\in K^l$ we set $W_{\bar{c}}(K) := \{ \bar{w} \in K^{4n}: (\bar{w},\bar{c})\in W(K) \}$. When $\bar{c}$ is not in the projection of $W(K)$ on $K^l$, $W_{\bar{c}}(K) = \emptyset$ and the conclusion of the theorem holds vacuously.} $\left(\bar{z},\bar{j}, \bar{j}', \bar{j}''\right) \in E^{\times}_{(z,J)}(K)\cap W_{\bar{c}}(K)$, if $$\dim W_{\bar{c}} < 3\dim S + \rk \Jac(\bar{z}),$$ then $\bar{j}\in T$ for some $T \in \Sigma$.
\end{theorem}

\begin{proof}
We may assume $S = (\Q^{\alg})^n$. First we show that there is a finite set $\Sigma(W,m)$ of special varieties, depending on $W$ and the number of derivations $m$, satisfying the conclusion of the theorem. Assume for some variety $W\subseteq (\Q^{\alg})^{4n+l}$ there is no such set $\Sigma(W,m)$. Consider the language $\{ +, \cdot, D_1, \ldots, D_m, \bar{z},\bar{j}, \bar{j}', \bar{j}'', \bar{c} \}$ where $\bar{z},\bar{j}, \bar{j}', \bar{j}'', \bar{c}$ are constant symbols. Let $\DF_0^m$ be the theory of differential fields of characteristic $0$ with $m$ commuting derivations. Consider the following set of sentence:
\begin{gather*}
    \DF_0^m \cup \{ D_i c_k =0: i, k \} \cup \{ (\bar{z},\bar{j}, \bar{j}', \bar{j}'')\in E^{\times}_{(z,J)} \cap W_{\bar{c}} \}\cup \\  \cup \{ \dim W_{\bar{c}} < 3n + \rk \Jac (\bar{z}) \} \cup \{ \Phi_N(j_i,j_k)\neq 0: i, k, N \}.
\end{gather*}
By our assumption every finite subset of the above set of sentences is satisfiable. Hence, by the compactness theorem of first-order logic, the whole set is satisfiable. This means that in some differential field the Ax--Schanuel theorem does not hold, which is a contradiction.

Now we show the existence of a set $\Sigma(W)$ with the desired properties which is independent of $m$. Let $(K; +, \cdot, D_1, \ldots, D_m)$ be a differential field with $m$ derivations and  $(\bar{z},\bar{j}, \bar{j}', \bar{j}'')\in E^{\times}_{(z,J)}(K) \cap W_{\bar{c}}(K)$. If $\rk \Jac(\bar{z}) = l$ then we can choose $l$ derivations, say $D_1,\ldots, D_l$, such that the matrix $(D_iz_k)_{1\leq i \leq l,1\leq k\leq n}$ has rank $l$. So we can work in the differential field $(K; +, \cdot, D_1, \ldots, D_l)$ which is a reduct of the original field. Since the number $l$ is bounded by $n$, we can take $\Sigma(W):= \bigcup_{1\leq l \leq n} \Sigma(W,l)$.
\end{proof}

As an immediate consequence we get the following result.

\begin{theorem}[Uniform Modular Ax--Schanuel without derivatives, cf. {\cite[Theorem 4.3]{Kirby-semiab}}]\label{uniform-Ax-no-der}
Let $S\seq (\Q^{\alg})^n$ be a $j$-special variety and let $W\seq (\Q^{\alg})^n \times S(\Q^{\alg}) \times (\Q^{\alg})^l$ be an algebraic variety defined over $\Q^{\alg}$. Then there is a finite  set $\Sigma = \Sigma(W)$ (depending on $W$) of proper $j$-special subvarieties of $S$ with the following property.

For every differential field $(K; D_1,\ldots,D_m)$ with $m$ commuting derivations and with field of constants $C$, for every  $\bar{c} \in C^l$ and for every $\left(\bar{z},\bar{j}\right) \in E^{\times}_{(z,j)}(K)\cap W_{\bar{c}}(K)$, if $\dim W_{\bar{c}} < \dim S + \rk \Jac(\bar{z})$, then $\bar{j}\in T$ for some $T \in \Sigma$.
\end{theorem}


\section{Weak Modular Zilber--Pink without derivatives}\label{weak-MZP}
\setcounter{equation}{0}

We begin by recalling the definition of (strongly) $j$-atypical subvarieties.

\begin{definition}
Let $V \subseteq S \subseteq C^n$ be an algebraic variety where $S$ is $j$-special. A \emph{$j$-atypical subvariety} of $V$ in $S$ is an irreducible component $W$ of some $V \cap T$, where $T$ is a $j$-special variety, such that $$\dim W> \dim V + \dim T -\dim S.$$
A $j$-atypical subvariety $W$ of $V$ is said to be \emph{strongly $j$-atypical} if no coordinate is constant on $W$.
\end{definition}

The following weak version of the Modular Zilber--Pink conjecture was proved by Pila and Tsimerman \cite[Theorem 7.1]{Pila-Tsim-Ax-j}. They use tools of o-minimality, while the proof that we give below is purely algebraic and is based on Kirby's adaptation of Zilber's proof of weak CIT (see \cite[Theorem 4.6]{Kirby-semiab}).

\begin{theorem}[Weak Modular Zilber--Pink]\label{weak-modularZP}
Given a parametric family of algebraic subvarieties $(V_{\bar{c}})_{\bar{c} \subseteq C}$ of a $j$-special variety $S$ in $C^n$, there is a finite collection $\Sigma$ of proper $j$-special subvarieties of $S$ such that for every $\bar{c} \subseteq C$, every strongly $j$-atypical subvariety of $V_{\bar{c}}$ in $S$ is contained in  some $T \in \Sigma$.
\end{theorem}

We will need the following concepts in the proof. 

\begin{definition}
The $j$-\textit{special closure} of an irreducible variety $X\subseteq C^n$ is the smallest $j$-special variety containing $X$.
\end{definition}

It is easy to see that irreducible components of an intersection of $j$-special varieties is $j$-special, hence there is a smallest $j$-special variety containing $X$.

\begin{definition}
A $C$-\emph{geodesic variety} $U\subseteq C^n$ (with coordinates $\bar{x}$) is an irreducible component of a variety defined by equations of the form $x_i=g_{i,k}x_k$ for some $g_{i,k}\in \SL_2(C)$. If $S\subseteq C^n$ (with coordinates $\bar{y}$) is a $j$-special variety, then $U$ is said to be a $C$-geodesic variety \emph{associated with} $S$ if for any $1\leq i, k\leq n$ we have $\Phi_N(y_i,y_k)=0$ on $S$ for some $N$ if and only if $x_i=g_{i,k}x_k$ on $U$ for some $g_{i,k}\in \SL_2(C)$.
\end{definition}

Note that for a $j$-special variety $S$ there are infinitely many geodesic varieties associated with $S$ since the matrices $g_{i,k}$ are chosen arbitrarily. Actually, the family of all geodesic varieties associated with $S$ forms a parametric family of varieties $(U_{\bar{c}})_{\bar{c} \subseteq C}$. In order to regard all geodesic varieties associated with all possible $j$-special varieties $T\subseteq C^n$ as members of a single parametric family of varieties we allow relations of the form $x_i=g_{i,k}x_k$ for $g_{i,k}=0$ (the zero matrix) which should be understood as the formula $0=0$ (i.e. we multiply through by a common denominator), that is, it does not impose any relations between $x_i$ and $x_k$. Thus, in a parametric family of geodesic varieties any two coordinates are related by an equation $x_i=g_{i,k}x_k$ where either $g_{i,k}\in \SL_2(C)$ or $g_{i,k}=0$.

Also, observe that if $U$ is a geodesic variety associated with a $j$-special $S$ then  $\dim U = \dim S$.

\begin{proof}[Proof of Theorem \ref{weak-modularZP}]
Let $W \subseteq V_{\bar{c}} \cap T$ be a strongly $j$-atypical subvariety of $V_{\bar{c}}$ where $T$ is a special subvariety of $S$. We know that
$$l:=\dim W > \dim V_{\bar{c}} + \dim T - \dim S.$$

Let $\bar{j} \in W$ be a generic point over $C$. We may assume that $T$ is the $j$-special closure of $\bar{j}$, i.e. the smallest $j$-special variety containing $\bar{j}$ (otherwise we would replace $T$ by the $j$-special closure of $\bar{j}$ and the above inequality would still hold). Consider the vector space $\Der(K/C)$ of derivations of the field $K:=C(\bar{j})$ over $C$. Its dimension is equal to $\td(K/C)$ which is equal to $\dim W$. Obviously, $\Der(K/C)$ is closed under the Lie bracket hence by Lemma \ref{Lie-Lemma} we can choose a commuting basis $D_1,\ldots, D_l$ of $\Der(K/C)$. Now let $U\subseteq K^n$ be a geodesic variety associated with $T$ which is defined by equations of the form $x_i=x_k$, i.e. all matrices $g\in \SL(C)$ occurring in the definition of $U$ are chosen to be the identity matrix (or the zero matrix if $x_i$ and $x_k$ are not linked). Pick a generic (over $K$) point $\bar{z}\in U$. Further, take a tuple $(\bar{j}', \bar{j}'')$ generic over $K(\bar{z})$ subject to the conditions that if $\Phi(j_i,j_k)=0$ for some modular polynomial $\Phi(Y_i,Y_k)$ then 
\begin{align*}
& \frac{\partial \Phi}{\partial Y_i}(j_i,j_k)\cdot j_i' + \frac{\partial \Phi}{\partial Y_k}(j_i,j_k)\cdot j_k' =0,\\
 & \begin{aligned}
 \frac{\partial^2 \Phi}{\partial Y_i^2}(j_i,j_k)\cdot \left(j_i'\right)^2+ \frac{\partial^2 \Phi}{\partial Y_k^2}(j_i,j_k)\cdot \left(j_k'\right)^2 + 2\cdot \frac{\partial^2 \Phi}{\partial Y_i \partial Y_k}(j_i,j_k)\cdot j_i'\cdot j_k'+ & \\ 
+\frac{\partial \Phi}{\partial Y_i}(j_i,j_k)\cdot j_i'' + \frac{\partial \Phi}{\partial Y_k}(j_i,j_k) & \cdot j_k''=  0. 
\end{aligned}
\end{align*}
These relations are obtained by differentiating the equation $\Phi(j_i,j_k)=0$ (see Proposition \ref{funct-eq}). Consider the field $F:=K(\bar{z},\bar{j}',\bar{j}'')$ and extend $D_i$ defining
$$D_iz_k = \frac{D_ij_k}{j_k'} = \frac{D_ij_k'}{j_k''} = \frac{D_ij_k''}{j_k'''},$$ where $j_k'''$ is uniquely determined from the equation $\Psi(j_k,j_k',j_k'',j_k''')=0$. By Proposition \ref{funct-eq} each $D_i$ is a derivation on $F$ and $(z_k,j_k)\in E_{(z,j)}(F)$ for all $k$. Straightforward calculations show that $D_i$'s commute on $F$ (cf. Claim \ref{Lambda-closed-Lie} in the proof of Theorem \ref{weak-MZPD}). Moreover, since no coordinate is constant on $W$, the elements $j_i$ are non-constant in the differential field $(F;+,\cdot, D_1, \ldots, D_l)$. It is also clear that $\rk \Jac(\bar{z}) = \rk \Jac(\bar{j}) = l$. Denote the field of constants of $F$ by $\hat{C}$.\footnote{Actually $C = \hat{C}$ but we do not need that in the proof.}


Observe that for different $j$-special varieties $T$ we get different geodesic varieties $U$. Nevertheless, all of those can be regarded as members of a parametric family of geodesic varieties $(U_{\bar{d}})_{\bar{d}\subseteq \hat{C} }$. Now if $U_{\bar{d}}$ is associated with $T$ then $\dim U_{\bar{d}} = \dim T$, and we have
$$\dim U_{\bar{d}} \times V_{\bar{c}} = \dim U_{\bar{d}} + \dim V_{\bar{c}} < \dim T+ \dim W +\dim S - \dim T = \rk \Jac (\bar{z}) + \dim S.$$

Now we apply the uniform Ax--Schanuel without derivatives to the parametric family $(U_{\bar{d}}\times V_{\bar{c}})_{\bar{c},\bar{d}\subseteq \hat{C}}$. There is a finite collection $\Sigma = \Sigma(V)$ of proper $j$-special subvarieties of $S$, depending on this parametric family only (which in turn depends only on the family $(V_{\bar{c}})_{\bar{c}\subseteq \hat{C}}$), such that $\bar{j}\in T$ for some $T\in \Sigma$. Since $\bar{j}$ is generic in $W$ over $C$ and $W$ is defined over $C$, and $T$ is defined over $\mathbb{Q}^{\alg}\subseteq C$, we must have $W\subseteq T$. 
\end{proof}

\begin{remark}
When the parametric family consists of a single variety $V$, we can choose the finite collection $\Sigma$ so that each $T \in \Sigma$ intersects $V$ strongly atypically in $S$. This can be achieved by repeatedly applying Theorem \ref{weak-modularZP}. It shows that $V$ contains only finitely many maximal strongly $j$-atypical subvarieties.
\end{remark}

The ideas of Section \ref{section-analytic-MZPD} can be used to give a complex analytic proof of Theorem \ref{weak-MZP}. It will be a mixture of the above proof and the proof of Pila and Tsimerman \cite[Theorem 7.1]{Pila-Tsim-Ax-j}, the main difference being that we do not use o-minimality and instead exploit the uniformity of differential Ax--Schanuel as in the above proof. We do not give more details as we prove a more general result using that method and it should be clear how this special case can be treated.

\section{D-special varieties}\label{D-special}
\setcounter{equation}{0}

\subsection{Definition and basic properties}
Now we define D-special varieties in $C^{3n}$ by analogy with $J$-special varieties (D stands for differential). The difference between D-special and $J$-special varieties is that we allow the geodesic relations to come from $\GL_2({C})$ rather than $\GL_2^+(\mathbb{Q})$.

\begin{definition}
Let $C$ be an algebraically closed field. Define $D$ as the zero derivation on $C$ and extend $(C;+,\cdot, D)$ to a differentially closed field $(K;+,\cdot, D)$. 
\begin{itemize}[leftmargin = 0.5cm]
    \item Let $T \subseteq C^n$ be a $j$-special variety and $U\subseteq C^n$ be a $C$-geodesic variety associated with $T$. Denote by $\langle \langle U, T \rangle \rangle$ the Zariski closure over $C$ of the projection of the set $$E^{\times}_{(z,J)}(K) \cap (U(K)\times T(K) \times K^2)$$ onto the last $3n$ coordinates.
    
    \item A D-\emph{special} variety is a variety $S:=\langle \langle U, T \rangle \rangle$ for some $T$ and $U$ as above. In this case $S$ is said to be a D-special variety associated with $T$ and $U$. We will also say that $T$ (or $U$) is a $j$-special (respectively, geodesic) variety associated with $S$. A D-special variety associated with $T$ is one associated with $T$ and $U$ for some $C$-geodesic variety $U$ which is associated with $T$. 
    
    \item $S \sim T$ means that $S$ is a D-special variety associated with $T$. For a set $\Sigma$ of $j$-special varieties $S \sim \Sigma$ means that $S\sim T$ for some $T \in \Sigma$.
    
    \item $\euscr{S}_{\D}$ is the collection of all D-special varieties.
\end{itemize}
\end{definition}

\begin{remarks}
\begin{itemize}[leftmargin = 0.5cm]

    \item[]

    \item Since the geodesic varieties associated with a fixed $T$ form a parametric family, we have a parametric family of D-special varieties associated with $T$.
    
    \item One can prove that $E^{\times}_{(z,J)}(K) \cap (U(K)\times T(K) \times K^2)$ is an irreducible Kolchin constructible set, which implies that $\langle \langle U, T \rangle \rangle$ is Zariski irreducible. Thus, D-special varieties are irreducible. We will actually prove this below by a different method. In Lemma \ref{lem: J_g special irreducible} we establish irreducibility of similar sets in an analytic context.
    
    \item It is clear that the definition does not depend on the ambient differentially closed field $K$. As we will see shortly, $\langle \langle U, T \rangle \rangle$ can be defined purely algebraically without referring to a differential field. 
    
    \item D-special varieties are automatically ``strongly'' special, i.e. they do not have any constant coordinates. In particular, a $j$-special variety associated with a D-special variety must be strongly $j$-special. In Section \ref{section-analytic-MZPD} we define more general special varieties over the complex numbers and show that they coincide with D-special varieties provided that they do not have constant coordinates. However, those special varieties are a generalisation of $J$-special varieties and may actually have constant coordinates. Nevertheless, in this paper we will only deal with strongly special varieties, hence any issues related to constant coordinates may be ignored.
\end{itemize}

\end{remarks}

As pointed out above, we can (and do) give an equivalent definition of D-special varieties using only algebraic language and completely avoiding mentioning derivations. Indeed, after choosing a geodesic variety $U$ associated with $T$, we can differentiate the modular relations defining $T$ and get some algebraic relations between $\bar{j},\bar{j}',\bar{j}''$, possibly over $\bar{z}$. Then eliminating $\bar{z}$ from those equations, that is, existentially quantifying over $\bar{z}$, we get algebraic equations defining $\langle \langle U, T \rangle \rangle$. 

Now let us discuss this strategy in more detail. 
We consider a simple case first. Let $T\subseteq C^2$ be a $j$-special variety defined by an equation $\Phi(y_1,y_2)=0$ where $\Phi$ is a modular polynomial. Let $U\subseteq C^2$ be a geodesic variety given by a single equation $x_2  = gx_1$ where $g=\begin{pmatrix} a & b \\ c & d \end{pmatrix} \in \SL_2(C)$. Now pick a non-constant point $$(z_1,z_2,j_1,j_2,j'_1,j'_2,j''_1,j''_2)\in K^8$$ such that $(z_i,j_i,j'_i, j''_i) \in E^{\times}_{(z,J)}(K)$ for $i=1,2$, and $\Phi(j_1,j_2)=0$ and $z_2=\frac{az_1+b}{cz_1+d}$. Applying\footnote{Recall that in a differential field $(K;+,\cdot, D)$ for a non-constant element $x\in K$ the derivation $\partial_x:K\rightarrow K$ is defined by $\partial_x: y \mapsto \frac{Dy}{Dx}$.} $\partial_{z_1}$ to the  equation $\Phi(j_1,j_2)=0$ we get
\begin{equation*}
    \frac{\partial \Phi}{\partial Y_1} (j_1,j_2) \cdot \partial_{z_1}j_1 +  \frac{\partial \Phi}{\partial Y_2} (j_1,j_2) \cdot \partial_{z_1}j_2 = 0.
\end{equation*}
On the other hand 
$$\partial_{z_1} = \frac{1}{(cz_1+d)^2}\cdot \partial_{z_2}.$$
Hence we have
\begin{equation}\label{eq-D-special}
    \frac{\partial \Phi}{\partial Y_1} (j_1,j_2) \cdot j'_1 + \frac{1}{(cz_1+d)^2}\cdot \frac{\partial \Phi}{\partial Y_2} (j_1,j_2) \cdot j'_2 = 0.
\end{equation}

Now if $c=0$, i.e. $g$ is upper triangular, then \eqref{eq-D-special} gives an algebraic relation between $j_1,j_2,j_1',j_2'$. Differentiating \eqref{eq-D-special} once more with respect to $z_1$ (i.e. applying $\partial_{z_1}$) we will get an algebraic relation between $j_1,j_2,j_1',j_2',j_1'',j_2''$ which, along with \eqref{eq-D-special} and the modular relation between $j_1$ and $j_2$, will define $\langle \langle U, T \rangle \rangle$. Indeed, the set defined by those equations is clearly irreducible, contains $\langle \langle U, T \rangle \rangle$ and has dimension $3$. By Ax--Schanuel, $\dim \langle \langle U, T \rangle \rangle \geq 3$, hence the above set is in fact equal to $\langle \langle U, T \rangle \rangle$.

If $c\neq 0$ then we do not get an algebraic relation between $j_1,j_2,j_1',j_2'$. Nevertheless we see that $z_1$ is algebraic over $j_1,j_2,j_1',j_2'$ and $j_2'$ is transcendental over $C(j_1,j_1',j_1'')$ since otherwise we would have $$\td_CC(z_1,j_1,j_1',j_1'')<4$$ which contradicts Ax--Schanuel. However, differentiating \eqref{eq-D-special} one more time we get an algebraic relation between $z_1,j_1,j_2,j_1',j_2',j_1'',j_2''$ which is linear with respect to $j_1''$ and $j_2''$. That equation, with \eqref{eq-D-special} and the equations defining $U$ and $T$, gives an irreducible subvariety of $C^8$. Therefore, its projection onto the last $6$ coordinates, as well as the Zariski closure of that, is an irreducible set of dimension $4$ which is equal to  $\langle \langle U, T \rangle \rangle$ as above.

Now assume that in addition to the above modular relation we also have a modular relation between $j_2$ and $j_3$ (we now work in $K^{12}$ and $U, T \subseteq K^3$). Note that this implies that $j_1$ and $j_3$ are also modularly dependent, and a modular equation between this coordinates is specified by $T$ as it is irreducible. Also, $z_1,z_2,z_3$ are pairwise linked by $\SL_2(C)$-relations. The above procedure can be used to write down the defining equations of $\langle \langle U, T \rangle \rangle$ in this setting. If all matrices from $\SL_2(C)$ linking $z_1,z_2,z_3$ are upper triangular then $\dim \langle \langle U, T \rangle \rangle = 3$, otherwise $\dim \langle \langle U, T \rangle \rangle =4$. Also, the equations are linear with respect to $j_i'$ or $j_i''$ for an appropriate $i$, which shows that $\langle \langle U, T \rangle \rangle$ is irreducible.

The same is also true for each subtuple of $\bar{j}$ the coordinates of which are pairwise modularly dependent; we can apply the above procedure to each such subtuple and get the equations defining $\langle \langle U, T \rangle \rangle$. For each such subtuple of maximal length $k$ we will get a distinct set of equations defining a subvariety of $K^{3k}$. We will refer to those as $j$-blocks.  

\begin{definition}
Let $V \subseteq C^{3n}$ be an irreducible variety. A \emph{$j$-block} of $V$ is a projection of $V$ onto the coordinates $(y_{i_1},\ldots,y_{i_k}, y'_{i_1},\ldots,y'_{i_k}, y''_{i_1},\ldots,y''_{i_k})$ for some $1\leq i_1<\ldots<i_k\leq n$ such that the coordinates $(y_{i_1},\ldots,y_{i_k})$ are pairwise modularly related on $V$ and none of them is modularly related to a coordinate $y_l$ for any $l\notin \{ i_1,\ldots,i_k\}$. The number $k$ is the size (or length) of the $j$-block. A $j$-block of length $1$ is called \textit{trivial}.
\end{definition}

This definition is more natural for D-special varieties but it is useful to have the concept of a $j$-block for arbitrary varieties. The above analysis shows that if $\langle \langle U, T \rangle \rangle$ consists of a single $j$-block then $\dim \langle \langle U, T \rangle \rangle = 3$ or $\dim \langle \langle U, T \rangle \rangle = 4$. Thus, we obtain the following characterisation of D-special varieties. 

\begin{proposition}
A D-special variety is irreducible and is equal to the product of its $j$-blocks, each of which has dimension $3$ (this corresponds to upper triangular matrices) or $4$. Furthermore, the number of the $j$-blocks of a D-special variety associated with a $j$-special variety $T$ is equal to $\dim T$.
\end{proposition}

Note also that a similar analysis for $J$-special varieties was carried out by Pila in his unpublished notes (see also \cite{Spence}). We will see in Section \ref{section-analytic-MZPD} that strongly $J$-special varieties are D-special.

\begin{definition}[cf. {\cite[Definition 1.3]{Spence}}]
Let $T\subseteq C^n$ be a $j$-special variety. A geodesic variety $U$ associated with $T$ is called \emph{upper triangular} if all matrices $g$ occurring in the definition of $U$ are upper triangular. If $U$ is upper triangular then a D-special variety associated with $T$ and $U$ is also called \emph{upper triangular}. 
\end{definition}

Observe that a D-special variety $S$ associated with a $j$-special variety $T$ is upper triangular if and only if $\dim S = 3 \dim T$. For example, $C^{3n}$ is an upper triangular D-special variety.

The following result often comes in useful.

\begin{lemma}\label{d-special-generic-J-point}
Assume $(K;+,\cdot,D)$ is an $\aleph_0$-saturated differentially closed field with field of constants $C$. Let $S$ be a D-special variety associated with a geodesic variety $U$, both defined over a finitely generated subfield $C_0 \seq C$. Then there is a point $(\bar{z},\bar{J})\in E^{\times}_{(z,J)}(K) \cap (U(K) \times S(K)) $ such that $\bar{J}$ is generic in $S$ over $C_0$ and $\bar{z}$ is generic in $U$ over $C_0$. In particular, $S$ contains a generic $E_J^{\times}(K)$-point over $C_0$.
\end{lemma}
\begin{proof}
The set $W:= E^{\times}_{(z,J)} \cap (U\times S)$ is Kolchin constructible, hence it can be decomposed into a union of irreducible relatively Kolchin closed subsets $W_i$.\footnote{As pointed out above, that set is actually Kolchin irreducible, so there is only one component. But it is not essential in the proof.} We may assume each $W_i$ is defined over $C_0$. The above observations on the structure of D-special varieties show that $\Zcl(W/C_0)$ is irreducible. On the other hand  $\Zcl(W/C_0) = \bigcup_i \Zcl(W_i/C_0) $, and each $\Zcl(W_i/C_0)$ is Zariski irreducible. Hence $\Zcl(W/C_0)$ must be equal to $\Zcl(W_i/C_0)$ for some $i$. Since $K$ is saturated and $W_i$ is Kolchin irreducible, $W_i(K)$ contains a Kolchin generic point $(\bar{z},\bar{J})\in E_{(z,J)}^{\times}(K)$ over $C_0$ and $\Zcl((\bar{z},\bar{J})/C_0) = \Zcl(W_i/C_0) = \Zcl(W/C_0)$. Therefore
$ \Zcl(\bar{J}/C_0) = S,~ \Zcl(\bar{z}/C_0) = U.$ 
\end{proof}

\begin{remark}
It can be proven that in the above lemma we can find a point generic not only over $C_0$, but also over $C$. This can be achieved, for example, using the methods of \cite{Aslanyan-adequate-predim} and \cite{Aslanyan-Eterovic-Kirby-Diff-EC-j}. Nevertheless, we do not prove it as we will not use it in this paper.
\end{remark}

\begin{proposition}\label{prop-S1=S2}
Let $S_1\subseteq S_2$ be D-special varieties associated with the same $j$-special variety $T$. Then $S_1 = S_2$.
\end{proposition}
\begin{proof}
We need to prove that each $j$-block of $S_1$ coincides with the corresponding $j$-block of $S_2$. Hence, we may assume both $S_1$ and $S_2$ have one $j$-block. If $\dim S_2 = 3$ then obviously $\dim S_1 = 3$ and $S_1=S_2$, so assume $\dim S_2 = 4$. By taking a non-upper triangular projection we may assume that $S_1 \subseteq S_2 \subseteq C^6$. In that case $T\subseteq C^2$ is defined by a single modular equation $\Phi(y_1,y_2)=0$. Also, assume $S_1$ and $S_2$ are defined over  $C_0\subseteq C$.

Extend $C$ to an $\aleph_0$-saturated differentially closed field $(K; +, \cdot, D)$. Let $\bar{J}:=(\bar{j},\bar{j}',\bar{j}'')\in S_1(K) \cap E_J^{\times}(K)$ be generic over $C_0$.

Differentiating the equality $\Phi(j_1,j_2)=0$ we get 
\begin{equation}\label{Dj_1 - Dj_2-0}
     \frac{\partial \Phi}{\partial Y_1} (j_1,j_2) \cdot Dj_1 + \frac{\partial \Phi}{\partial Y_2} (j_1,j_2) \cdot Dj_2 = 0.
\end{equation}

Further, we have
\begin{equation}\label{Dj/j'-0}
    Dj_1' = j_1''\cdot \frac{Dj_1}{j_1'} ~ \mbox{  and  }~ Dj_2' = j_2''\cdot \frac{Dj_2}{j_2'}.
\end{equation}

Let $U_2\subseteq C^2$ be a $C$-geodesic variety associated with $S_2$, defined by an equation
$$x_2 = \frac{ax_1+b}{cx_1+d}, \mbox{ where } a, b, c, d \in C,~ ad-bc =1,~ c \neq 0.$$ Pick $z_1 \in K$ such that

\begin{equation}\label{eq-D-special-2-0}
    \frac{\partial \Phi}{\partial Y_1} (j_1,j_2) \cdot j'_1 + \frac{1}{(cz_1+d)^2}\cdot \frac{\partial \Phi}{\partial Y_2} (j_1,j_2) \cdot j'_2 = 0.
\end{equation}

Differentiating this equality and using \eqref{Dj_1 - Dj_2-0} and \eqref{Dj/j'-0}, we obtain
\begin{equation*}
    Dz_1 = \xi(z_1, j_1, j_2, j_1', j_2', j_1'', j_2'')\cdot Dj_1
\end{equation*}
where $\xi$ is a rational function. We claim that $$ \xi(z_1, j_1, j_2, j_1', j_2', j_1'', j_2'') = \frac{1}{j_1'}.$$

Let $$(\bar{u},\bar{v},\bar{v}',\bar{v}'') \in E_{(z,J)}(K) \cap (U_2(K)\times S_2(K))$$ be generic over $C_0$ as in Lemma \ref{d-special-generic-J-point}. Differentiating $\Phi(v_1,v_2)=0$ with respect to $D$ and $\partial_{u_1}$ we get 

\begin{equation*}
     \frac{\partial \Phi}{\partial Y_1} (v_1,v_2) \cdot D v_1 + \frac{\partial \Phi}{\partial Y_2} (v_1,v_2) \cdot D v_2 = 0
\end{equation*}
and
\begin{equation*}
    \frac{\partial \Phi}{\partial Y_1} (v_1,v_2) \cdot v'_1 + \frac{1}{(cu_1+d)^2}\cdot \frac{\partial \Phi}{\partial Y_2} (v_1,v_2) \cdot v'_2 = 0.
\end{equation*}

Differentiating the second equality and taking into account the fact that $$\frac{D v_1}{v_1'} = \frac{D v_1'}{v_1''} \mbox{ and } \frac{D v_2}{v_2'} = \frac{D v_2'}{v_2''}$$ we see that
\begin{equation*}
    \delta u_1 = \xi(u_1, v_1, v_2, v_1', v_2', v_1'', v_2'')\cdot \delta v_1
\end{equation*}
where $\xi$ is the same rational function as above. However, we know that $\delta u_1 = \frac{\delta v_1}{v'_1}$ and hence $\xi(u_1, v_1, v_2, v_1', v_2', v_1'', v_2'') = \frac{1}{v'_1}$. Since $(v_1, v'_1, v''_1,v_2, v'_2, v''_2)$ is generic in $S_2$ over $C_0$ (and $u_1$ and $z_1$ satisfy the same algebraic equation over $v_1,v_2,v_1',v_2'$ and $j_1,j_2,j_1',j_2'$ respectively), we conclude that $\xi(z_1, j_1, j_2, j_1', j_2', j_1'', j_2'') = \frac{1}{j'_1}$, that is, 
$$Dz_1 = \frac{Dj_1}{j_1'}.$$ 

Now if $z_2:=\frac{az_1+b}{cz_1+d}$ then it is clear that $j_i = \partial_{z_i}j_i,~ j_i'' = \partial_{z_i}j_i',~ j_i''' = \partial_{z_i}j_i''$ and $(z_i,j_i,j_i',j_i'')\in E_{(z,J)}(K)$ for $i=1,2$. Therefore 
$$\dim S_1 \geq \td_CC(\bar{J}) = \td_CC(\bar{z},\bar{J})\geq 4.$$
Thus $\dim S_1=\dim S_2$ and $S_1\subseteq S_2$, hence $S_1=S_2$.
\end{proof}

\begin{corollary}\label{cor-D-special-EJ-generic}
Let $K$ be a differential field, $S$ be a D-special variety defined over the field of constants and $\bar{J}\in S(K)\cap E^{\times}_J(K)$. Then the appropriate projections of $\bar{J}$ are generic in $j$-blocks of $S$.
\end{corollary}

\subsection{D-special closure}

\begin{definition}
Let $V\subseteq C^{3n}$ be an algebraic variety (or, more generally, an arbitrary set). \emph{A D-special closure} of $V$ is a D-special variety $S\subseteq C^{3n}$ which contains $V$ and is minimal among the D-special varieties containing $V$. 
\end{definition}

\begin{remark}
By Noetherianity of the Zariski topology every variety has at least one D-special closure which, in general, is not unique. When all $j$-blocks of a variety $V$ are D-special, $V$ has a unique D-special closure which is the product of its $j$-blocks and is the smallest D-special variety containing $V$. 
\end{remark}

\begin{example}
Let $T_1, T_2 \subseteq C^6$ be D-special varieties of dimension $4$ (with a single $j$-block). Further, let $S_1 := C^3 \times T_1 \subseteq C^9$ be a D-special variety whose projection on $(y_1,y_1',y_1'')$ is $C^3$ and the projection on $(y_2,y_3,y'_2,y'_3,y''_2,y''_3)$ is equal to $T_1$. Similarly defined $S_2 := T_2 \times C^3 \subseteq C^9$ whose projection on $(y_3,y_3',y_3'')$ is $C^3$ and the projection on $(y_1,y_2,y'_1,y'_2,y''_1,y''_2)$ is equal to $T_2$. If $V\subseteq S_1\cap S_2$ is a component of the intersection, then $\dim W = 5$, so it cannot be D-special. Both $S_1$ and $S_2$ (each of dimension $7$) are D-special closures of $W$.
\end{example}

\begin{lemma}\label{lemma-over-C}
Let $(K;+, \cdot, D)$ be a differential field with an algebraically closed field of constants $C$, and let $\bar{J}=(\bar{j},\bar{j}', \bar{j}'')\in E^{\times}_J(K)$. If $T\subseteq K^n$ is the $j$-special closure of $\bar{j}$ then $\bar{J}$ belongs to a D-special variety $S \sim T$ defined over $C$.
\end{lemma}
\begin{proof}
Let $F$ be a differential closure of $K$. Since $C$ is algebraically closed, it is the field of constants of $F$. Pick $\bar{z}\in F^n$ such that $(\bar{z},\bar{J})\in E^{\times}_{(z,J)}(F)$. Then $\bar{z}$ belongs to a $C$-geodesic variety associated with $T$, hence $\bar{J}$ belongs to a D-special variety defined over $C$ and associated with $T$.
\end{proof}

\begin{corollary}
Let $(K;+, \cdot, D)$ be a differential field with an algebraically closed field of constants $C$ and $\bar{J}=(\bar{j},\bar{j}', \bar{j}'')\in E^{\times}_J(K)$. Then $\bar{J}$ has a unique D-special closure and it does not depend on $K$, i.e. if $F$ is a differential superfield of $K$ with field of constants $\hat{C}$ then the D-special closure of $\bar{J}$ over $\hat{C}$ is equal to that over $C$.
\end{corollary}

\subsection{Weak Ax--Schanuel for D-special varieties}

The following is a weak form of Ax--Schanuel in terms of D-special varieties.

\begin{theorem}[Weak Ax--Schanuel]\label{weak-Ax}
Assume $(K;+,\cdot, D_1,\ldots, D_m)$ is a differential field (with commuting derivations) with constant field $C$. Let $S \subseteq C^{3n}$ be a D-special variety associated with a $j$-special variety $T\subseteq K^n$. If $\bar{J}:=(\bar{j}, \bar{j}', \bar{j}'') \in E^{\times}_J(K) \cap S(K)$ and $$\td_CC(\bar{J})< \dim S - \dim T + \rk \Jac(\bar{j})$$ then $\bar{j}$ belongs to a proper $j$-special subvariety of $T$ and hence  $\bar{J}$ belongs to a proper D-special subvariety of $S$.
\end{theorem}
\begin{proof}
Assume $\bar{j}$ does not belong to a proper $j$-special subvariety of $T$, that is, $T$ is the $j$-special closure of $\bar{j}$. By Proposition \ref{prop-S1=S2}, $S$ is the D-special closure of $\bar{J}$. 

Pick a tuple $\bar{z}$ (possibly in a differential field extension $F$ of $K$) such that $(\bar{z},\bar{J})\in E^{\times}_{(z,J)}(F)$. It is clear that if $t$ is the number of $j$-blocks of $S$ of dimension $3$ then $$\td(C(\bar{z},\bar{J})/C(\bar{J})) \leq t $$ and $$\td_CC(\bar{z},\bar{J}) < t + \dim S - \dim T + \rk \Jac(\bar{j}) = 3 \dim T + \rk \Jac(\bar{j}).$$ This contradicts the Ax--Schanuel theorem.
\end{proof}

A uniform version of this theorem can also be proved as before.

\section{An analytic approach}\label{section-analytic-MZPD}
\setcounter{equation}{0}

\subsection{$J_{\bar{g}}$-special varieties}

For $g\in \GL_2(\mathbb{C})$ let $\mathbb{H}^g := g^{-1}\mathbb{H}$ and let $j_g: \mathbb{H}^g \rightarrow \mathbb{C}$ be the function $j_g(z) = j(gz)$. For a tuple $\bar{g} = (g_1, \ldots, g_n) \in \GL_2(\mathbb{C})^n$ set $\mathbb{H}^{\bar{g}}:= \mathbb{H}^{g_1}\times \cdots \times \mathbb{H}^{g_n}$ and define functions
$$j_{\bar{g}}: \mathbb{H}^{\bar{g}}\rightarrow \mathbb{C}^n : (z_1,\ldots,z_n) \mapsto (j_{g_1}(z_1), \ldots, j_{g_n}(z_n))$$
and $$J_{\bar{g}} = (j_{\bar{g}},j'_{\bar{g}},j''_{\bar{g}}): \mathbb{H}^{\bar{g}}\rightarrow \mathbb{C}^{3n}: \bar{z} \mapsto (j_{\bar{g}}(\bar{z}),j'_{\bar{g}}(\bar{z}),j''_{\bar{g}}(\bar{z}))$$ where the derivation is coordinatewise and $$j'_{g_i}(z_i) =\frac{\dd}{\dd z_i} j(g_iz_i) = j'(g_iz_i)\cdot g_i'z_i,~ j_{g_i}''(z_i) = \frac{\dd^2}{\dd z_i^2}j(g_iz_i).$$ Here $j'(z) = \frac{\dd}{\dd z}j(z)$ and $g'z = \frac{\dd}{\dd z}(gz)$ for $g\in \GL_2(\mathbb{C})$. We let $\Gamma_{\bar{g}}\subseteq \mathbb{H}^{\bar{g}}\times \mathbb{C}^{3n}$ be the graph of $J_{\bar{g}}$. When $g_i$ is the indetity matrix for each $i$, we drop the subscript $\bar{g}$ from $J_{\bar{g}}$ and $\Gamma_{\bar{g}}$.

\begin{definition}
Let $\bar{g} = (g_1, \ldots, g_n) \in \GL_2(\mathbb{C})^n$. An $\mathbb{H}^{\bar{g}}$-\emph{special} variety is an irreducible component of a subvariety of $\mathbb{H}^{\bar{g}}$ defined by equations of the form
\begin{equation}\label{eq-Jg-special}
    z_i = g_i^{-1} \gamma_{i,k} g_k z_k
\end{equation}
where $\gamma_{i,k}\in \GL_2^+(\mathbb{Q})$.  
\end{definition}

\begin{definition}
Let $\bar{g}\in \GL_2(\mathbb{C})^n$.
\begin{itemize}[leftmargin = 0.5cm]
    \item For an $\mathbb{H}^{\bar{g}}$-special variety $U \subseteq \mathbb{H}^{\bar{g}}$ denote the set $\Zcl(J_{\bar{g}}(U)/\mathbb{C})$ by $\langle \langle U \rangle \rangle_{\bar{g}}$.
    \item A $J_{\bar{g}}$-\emph{special} variety is a variety of the form $\langle \langle U \rangle \rangle_{\bar{g}}$ for some $\mathbb{H}^{\bar{g}}$-special $U$.
    \item A $J_{\bar{g}}$-special variety $S$ is \emph{strongly} $J_{\bar{g}}$-special if no coordinate is constant on $S$.
    \item For a $J_{\bar{g}}$-special variety $S$, we say it is associated with a $j$-special variety $T$, and write $S \sim T$, if $S= \langle \langle U \rangle \rangle_{\bar{g}}$ for some $U\subseteq \mathbb{H}^{\bar{g}}$ with $j_{\bar{g}}(U) = T$ (equivalently, the projection of $S$ onto the first $n$ coordinates is equal to $T$).
    \item For a set $\Sigma$ of $j$-special varieties $S\sim \Sigma$ means that $S\sim T$ for some $T\in \Sigma$.
    \item The collection of all strongly $J_{\bar{g}}$-special varieties is denoted by $\euscr{S}_{\bar{g}}$.
\end{itemize}
\end{definition}

\begin{remarks}
\begin{itemize}[leftmargin = 0.5cm]
    \item []
    \item A variety can simultaneously be $J_{\bar{g}}$-special and $J_{\bar{h}}$-special for some $\bar{h}\neq \bar{g}$, for only the coordinates occurring in \eqref{eq-Jg-special} are relevant. If a $J_{\bar{g}}$-special variety $S$ contains a $J_{\bar{h}}$-special variety, then $S$ is also $J_{\bar{h}}$-special.

    \item It is clear that when $g_i$ is the identity matrix for each $i$ then strongly $J_{\bar{g}}$-special varieties coincide with strongly $J$-special varieties. However, $\bar{J}$-special varieties are defined over $\mathbb{Q}^{\alg}$, while $J_{\bar{g}}$-special varieties are defined over $\mathbb{C}$, so in general these two notions are not the same. This distinction does not cause any issues because we work only with strongly special varieties.
\end{itemize}
\end{remarks}

\begin{lemma}\label{lem: J_g special irreducible}
$J_{\bar{g}}$-special varieties are irreducible.
\end{lemma}
\begin{proof}
Let $S = \langle \langle U \rangle \rangle_{\bar{g}}$ be $J_{\bar{g}}$-special. By definition, $U$ is an irreducible complex analytic set, hence its image $J_{\bar{g}}(U)$ is also irreducible (see the footnote to Remark \ref{remark-image-analytic}). If $\langle \langle U \rangle \rangle_{\bar{g}}$ is reducible, then one of its algebraically irreducible components must contain $J_{\bar{g}}(U)$, therefore that component is equal to the Zariski closure of the latter. Thus, $\langle \langle U \rangle \rangle_{\bar{g}}$ is equal to an irreducible component of itself, hence it is in fact irreducible.
\end{proof}

\begin{proposition}
A subvariety of $\mathbb{C}^{3n}$ is D-special if and only if it is strongly $J_{\bar{g}}$-special for some $\bar{g}\in \GL_2(\mathbb{C})^n$. In other words,  $\euscr{S}_{\D}  = \bigcup_{\bar{g}\in \GL_2(\mathbb{C})^n}\euscr{S}_{\bar{g}}.$
\end{proposition}
\begin{proof}
Let $U$ be $\mathbb{H}^{\bar{g}}$-special defined by equations \eqref{eq-Jg-special}. It is clear that $T:=j_{\bar{g}}(U)$ is a $j$-special variety defined by equations $\Phi_N(y_i, y_k)=0$ with $N = N(\gamma_{i,k})$ (see Section \ref{jBackground}). 
Let $W \subseteq \mathbb{C}^n$ be the $\mathbb{C}$-geodesic variety (associated with $T$) defined by the same equations as $U$, that is, $U = W \cap \mathbb{H}^{\bar{g}}$. Then it is straightforward to show that $\langle \langle U \rangle \rangle_{\bar{g}} = \langle \langle W, T \rangle \rangle$.

For the converse, assume $S\subseteq \mathbb{C}^{3n}$ is a D-special variety associated with $T$ and $W$ defined over a finitely generated subfield $C_0 \subseteq \mathbb{C}$. Let $K\supseteq \mathbb{C}$ be an $\aleph_0$-saturated differentially closed field with constant field $\mathbb{C}$. By Lemma \ref{d-special-generic-J-point} there is a point $(\bar{z},\bar{J})\in E^{\times}_{(z,J)}(K) \cap (W(K) \times S(K)) $ such that $\bar{J}$ is generic in $S$ over $C_0$ and $\bar{z}$ is generic in $W$ over $C_0$. Let $K_0 \subseteq K$ be a finitely generated differential field, with constant field $C_0$, containing $\bar{J}$ and $\bar{z}$. By Seidenberg's embedding theorem $K_0$ can be embedded into a differential field of meromorphic functions on a complex domain. Then by Lemma \ref{solutions-of-eq-j} $j_i = j(g_iz_i)$ where $j$ is the $j$-function and $g_i \in \GL_2(\mathbb{C})$. In other words, $j_i = j_{g_i}(z_i)$ and $\bar{J} = J_{\bar{g}}(\bar{z})$ where $\bar{g}=(g_1,\ldots,g_n)$. Let $U:= \mathbb{H}^{\bar{g}}\cap W$. Then it is easy to verify that $U$ is $\mathbb{H}^{\bar{g}}$-special and $S=\langle \langle U \rangle \rangle_{\bar{g}}$.
\end{proof}

This shows, in particular, that the structure of $J_{\bar{g}}$-special varieties is similar to that of D-special varieties, that is, a $J_{\bar{g}}$-special variety is equal to the product of its $j$-blocks. However, since $J_{\bar{g}}$-special varieties may have constant coordinates, their $j$-blocks may be of dimension zero. Of course, the $j$-blocks of a strongly $J_{\bar{g}}$-special variety are of dimension $3$ or $4$ depending on the matrices involved in the defining equations of the corresponding $\mathbb{H}^{\bar{g}}$-special variety.


\begin{notation}
Let $\pr_j : \mathbb{C}^{3n} \rightarrow \mathbb{C}^n$ be the projection onto the $j$-coordinates, i.e. the first $n$ coordinates. Similarly, let $\ppr_j : \mathbb{C}^{4n} \rightarrow \mathbb{C}^n$ be the projection onto the second $n$ coordinates.
\end{notation}

The following is an equivalent form of the Complex Ax--Schanuel for the function $J$ \cite[Theorem 1.2]{Pila-Tsim-Ax-j}. 

\begin{theorem}\label{complex-Ax-Sch}
Let $V\subseteq \mathbb{C}^{4n}$ be an algebraic variety and let $A$ be an analytic component of the intersection $V \cap \Gamma$. If $\dim A > \dim V - 3n$ and no coordinate is constant on $\ppr_j A$ then it is contained in a proper $j$-special subvariety of $\mathbb{C}^n$.
\end{theorem}

We will need the following uniform version of this theorem.

\begin{theorem}\label{complex-Ax-Sch-uniform}
Let $S\subseteq \mathbb{C}^{3n}$ be an upper triangular D-special variety, associated with a $j$-special $T$, and $V_{\bar{c}}\subseteq \mathbb{C}^n\times S(\mathbb{C})$ be a parametric family of algebraic varieties. Then there is a finite collection $\Sigma$ of proper $j$-special subvarieties of $T$ such that for every $\bar{c}\subseteq \mathbb{C}$ and every $\bar{g}\in \GL_2(\mathbb{C})^n$, if $A$ is an analytic component of the intersection $V_{\bar{c}} \cap \Gamma_{\bar{g}}$ with $\dim A > \dim V_{\bar{c}} - \dim S$, and no coordinate is constant on $\ppr_j A$, then $\ppr_j A$ is contained in some $T' \in \Sigma$.
\end{theorem}
\begin{proof}
We replicate the argument of \cite[$\S 2.5$]{Pila-Tsim-Ax-j}. Let $(\bar{z}, j_{\bar{g}}(\bar{z}), j_{\bar{g}}'(\bar{z}), j_{\bar{g}}''(\bar{z}))$ be local coordinates on $A$ where each $z_i = z_i(w_1, \ldots, w_l)$ is a holomorphic function of $\bar{w}$ defined on some open subset $W$ of $\mathbb{C}^l$. Consider the field $K$ of meromorphic functions on $W$, equipped with derivations $D_i = \frac{\dd}{\dd w_i},~ i=1, \ldots, l$. Then $(\bar{z}, j_{\bar{g}}(\bar{z}), j_{\bar{g}}'(\bar{z}), j_{\bar{g}}''(\bar{z}))\in E_{(z,J)}^{\times}(K)$, $\dim A = \rk (D_i z_k)_{i,k}$ and $\dim V_{\bar{c}} < \dim S + \rk (D_i z_k)_{i,k} = 3 \dim T + \rk (D_i z_k)_{i,k}$. By Theorem \ref{uniform-AS-j} there is a finite set $\Sigma$ of proper $j$-special subvarieties of $T$, depending on $V$ but not on $\bar{c}$ (neither on $\bar{g}$), such that $j_{\bar{g}}(\bar{z})\in T'$ for some $T' \in \Sigma$. But then $\ppr_j A \subseteq T'$.
\end{proof}

\subsection{Weak Zilber--Pink for $J_{\bar{g}}$}

\begin{definition}\label{def-an-str-atyp}
For $J_{\bar{g}}$-special varieties $T \subseteq S\subseteq \mathbb{C}^{3n}$ and a subvariety $V \subseteq S$ an atypical (in $S$) component $X$ of the intersection $V \cap T$ is  \emph{strongly $J_{\bar{g}}$-atypical} if for every irreducible analytic component $Y$ of $X \cap J_{\bar{g}}(\mathbb{H}^{\bar{g}})$ no coordinate is constant on the projection $\pr_j Y$. The \emph{strongly $J_{\bar{g}}$-atypical set} of $V$ in $S$, denoted $\SAtyp_{\bar{g}}(V;S)$, is the union of all strongly $J_{\bar{g}}$-atypical subvarieties of $V$. 
\end{definition}

\begin{theorem}\label{weak-MZPD1}
Let $S\subseteq \mathbb{C}^{3n}$ be an upper triangular D-special variety associated with a $j$-special variety $T\subseteq \mathbb{C}^n$. For a parametric family of algebraic varieties $V_{\bar{c}} \subseteq S$ there is a finite collection $\Sigma$ of proper $j$-special subvarieties of $T$ such that for every $\bar{c}$ and every $\bar{g}\in \GL_2(\mathbb{C})$
$$\SAtyp_{\bar{g}}(V_{\bar{c}};S) \cap J_{\bar{g}}(\mathbb{H}^{\bar{g}}) \subseteq \bigcup_{\substack{P\sim \Sigma\\ P \in \euscr{S}_{\bar{g}}}} P.$$ 
\end{theorem}

This implies Theorem \ref{weak-MZPD1-intro}, and itself follows from the next proposition.

\begin{proposition}\label{prop-WMZPD1}
Let $S\subseteq \mathbb{C}^{3n}$ be an upper triangular D-special variety associated with a $j$-special variety $T\subseteq \mathbb{C}^n$. For a parametric family of algebraic varieties $V_{\bar{c}} \subseteq S$ there is a finite collection $\Sigma$ of proper $j$-special subvarieties of $T$ such that for every $\bar{c}$, every $\bar{g}\in \GL_2(\mathbb{C})$ and every $J_{\bar{g}}$-atypical subvariety $X$ of $V_{\bar{c}}$, if $A \subseteq X \cap J_{\bar{g}}(\mathbb{H}^{\bar{g}})$ is an analytic component with $\pr_j A$ having no constant coordinates, then $\pr_j A$ is contained in some $T' \in  \Sigma$.
\end{proposition}

We will need some auxiliary results in the proof of this proposition. First, we will need the following theorem on the dimension of intersection in analytic sets.

\begin{theorem}\label{dimension-intersect-nonsmooth}
Let $A,B\subseteq M$ be irreducible analytic varieties, and $X$ be an irreducible component of $A \cap B$. If $X$ contains a non-singular point of $M$ then $$\dim X \geq \dim A + \dim B - \dim M.$$
\end{theorem}
\begin{proof}
This is a well-known theorem for smooth $M$ (see \cite[Chapter III, 4.6]{Loj-analytic-geom}). If $M$ is not smooth, then we can intersect all varieties with the set of non-singular points of $M$, which will not change their dimensions (in particular, those intersections will be non-empty due to the hypothesis of the theorem).
\end{proof}

In order to apply this theorem, we need the following result.

\begin{lemma}\label{10-lemma-smooth}
Assume $T\subseteq \mathbb{C}^{3n}$ is $J_{\bar{g}}$-special and $Y$ is a complex analytically irreducible subset of $T \cap J_{\bar{g}}(\mathbb{H}^{\bar{g}})$ such that no coordinate is constant on the projection $\pr_j Y$. Then $Y$ contains a non-singular point of $T$.
\end{lemma}
\begin{proof}
Let $T_{\s} \subseteq T$ be the subset of singular points of $T$. Then $T_{\s}$ is a proper Zariski closed subset of $T$. We need to show that $Y \nsubseteq T_{\s}$. 

Let us assume first that $T$ consists of a single $j$-block. We claim that all but countably many points of $Y$ are non-singular points of $T$. Consider the set $Z:= J_{\bar{g}}^{-1} (Y\cap T_{\s})$. If $Z$ is at most countable, then we are done. Otherwise for some coordinate $z_i$ the projection of $Z$ on $z_i$, denoted $Z_i$, must be uncountable. Therefore, $Z_i$ contains a limit point of itself. On the other hand, the functions $z_i, j_{g_i}(z_i), j_{g_i}'(z_i), j_{g_i}''(z_i)$ must satisfy a non-trivial algebraic equation for $z_i \in Z_i$ as $T_{\s} \subsetneq T$. Since $Z_i$ has a limit point, that identity will hold for all $z_i\in \mathbb{H}^{g_i}$ which contradicts Mahler's theorem (Ax--Schanuel for $n=1$). 

Now let $T= T_1 \times \ldots \times T_k$ be the decomposition of $T$ into a product of $j$-blocks. By the above argument, the projection of $Y$ on $T_i$ contains at most countably many singular points of $T_i$. Denote that set by $S_i$. If a point is smooth on each $T_i$ then it will be smooth on $T$ as well. So assume, for contradiction, that every point of $Y$ is singular on at least one $T_i$. It means that for every point of $Y$, its projection on some $T_i$ is contained in $S_i$. But then $Y$ will be contained in a countable union of varieties each of which has a constant coordinate. Hence $Y$ must be contained in one of them, which implies that the projection $\pr_j Y$ has a constant coordinate, and this is a contradiction.
\end{proof}

\begin{lemma}\label{lemma-component-J(U)}
Let $S$ be a $J_{\bar{g}}$-special variety and $Y$ be an analytic component of $S \cap J_{\bar{g}}(\mathbb{H}^{\bar{g}})$. Then there is an $\mathbb{H}^{\bar{g}}$-special variety $U$ such that $Y = J_{\bar{g}}(U)$ and  $S = \langle \langle U \rangle \rangle_{\bar{g}}$.
\end{lemma}
\begin{proof}
Let $(j_{\bar{g}}(\bar{z}), j_{\bar{g}}'(\bar{z}), j_{\bar{g}}''(\bar{z}))$ be local coordinates on $Y$ where each $z_i = z_i(\bar{w})$ is a function of $\bar{w}$. Clearly, there is an $\mathbb{H}^{\bar{g}}$-special variety $U$ such that $j_{\bar{g}}(U)= \pr_j S$ and $\bar{z}(\bar{w})\in U$ for all $\bar{w}$. Since $Y$ is irreducible and $Y\subseteq S$, we have $Y \subseteq J_{\bar{g}}(U) \subseteq S$. Therefore $Y = J_{\bar{g}}(U)$, for $Y$ is a component of $S \cap J_{\bar{g}}(\mathbb{H}^{\bar{g}})$. It is also clear that $S = \langle \langle U \rangle \rangle_{\bar{g}}$.
\end{proof}


\begin{proof}[Proof of Proposition \ref{prop-WMZPD1}]
We consider the case of a single variety $V$ first. 
Let $T \subseteq S$ be a $J_{\bar{g}}$-special variety and $X \subseteq V \cap T$ be an atypical component in $S$. Assume $A \subseteq X \cap J_{\bar{g}}(\mathbb{H}^{\bar{g}})$ is an analytic component such that no coordinate is constant on $\pr_j A$. Since $A \subseteq T \cap J_{\bar{g}}(\mathbb{H}^{\bar{g}})$, by Lemma \ref{lemma-component-J(U)} there is an $\mathbb{H}^{\bar{g}}$-special variety $U$ such that $A \subseteq J_{\bar{g}}(U)\subseteq T$. Thus, $A$ is an analytic component of  $X \cap J_{\bar{g}}(U)$. 
By Lemma \ref{10-lemma-smooth} $A$ contains a non-singular point of $T$. So by Theorem \ref{dimension-intersect-nonsmooth} we have\footnote{As pointed out earlier, the set $J_{\bar{g}}(U)$ may not be an analytic set. However, it is locally analytic, that is, every point of $J_{\bar{g}}(U)$ has a Euclidean neighbourhood where the set is analytic. This allows one to use the theorem on dimension of intersection of analytic sets. Moreover, if we restrict $j$ to the interior of a fundamental domain, then $j$ becomes a bijective holomorphic map, hence has a holomorphic inverse. Then the intersection of $J_{\bar{g}}(U)$ with a product of such fundamental domains will be a globally analytic subset of the latter, whose dimension is equal to $\dim J_{\bar{g}}(U)$. This also shows that we can use Theorem \ref{dimension-intersect-nonsmooth}.} 
\begin{gather*}
    \dim A  \geq \dim X + \dim J_{\bar{g}}(U) - \dim T >\\
   \dim V + \dim T - \dim S + \dim J_{\bar{g}}(U) - \dim T = \dim V + \dim U - \dim S.
\end{gather*}
This implies
$$ \dim ((U\times A) \cap \Gamma_{\bar{g}}) = \dim A > \dim (U \times V) - \dim S.$$ 
Now the desired result follows from Theorem \ref{complex-Ax-Sch-uniform} applied to the parametric family of algebraic varieties $W_{\bar{c}} \times V$ where $W_{\bar{c}}$ varies over the parametric family of all $\mathbb{C}$-geodesic varieties.

It is clear that the proof goes through for a parametric family $V_{\bar{c}}$.
\end{proof}

Since strongly $J$-special subvarieties of an algebraic variety $V$ are obviously strongly $J$-atypical, Theorem \ref{weak-MZPD1} implies the following weak version of the MAOD conjecture, which also follows from Spence's results \cite{Spence}. In Section \ref{weak-MAOD} we state and prove a more general Functional MAOD result.

\begin{theorem}
For every algebraic variety $V \subsetneq \mathbb{C}^{3n}$ there is a finite collection $\Sigma$ of proper $\mathbb{H}$-special subvarieties of $\mathbb{H}^n$ such that every strongly $J$-special subvariety of $V$ is contained in a $J$-special variety of the form $\langle \langle \bar{\gamma} U \rangle \rangle$ for some $\bar{\gamma}\in \SL_2(\mathbb{Z})^n$ and some $U \in \Sigma$.
\end{theorem}

\section{Differential Modular Zilber--Pink with Derivatives}\label{section-DMZPD}

\begin{definition}
For a D-special variety $S\subseteq C^{3n}$ and a subvariety $V \subseteq S$ we let the D-\emph{atypical set} of $V$ in $S$, denoted $\Atyp_{\D}(V;S)$, be the union of all D-atypical subvarieties of $V$ in $S$, that is, atypical components of intersections $V \cap T$ in $S$ where $T\subseteq S$ is D-special.  
\end{definition}

Now we formulate and prove a differential analogue of the MZPD conjecture, which is a general version of Theorem \ref{DMZPD-intro}.

\begin{theorem}[DMZPD]\label{DMZPD}
Let $(F; +,\cdot, D)$ be a differential field with an algebraically closed field of constants $C$. Let also $S\subseteq C^{3n}$ be an upper triangular D-special variety defined over $C$ associated with a $j$-special variety $T\subseteq C^n$. Given a parametric family of algebraic subvarieties $(V_{\bar{c}})_{\bar{c} \subseteq C}$ of $S$, there is a finite collection $\Sigma$ of proper $j$-special subvarieties of $T$ such that for every $\bar{c} \subseteq C$ we have
$$\Atyp_{\D}(V_{\bar{c}};S)(F) \cap E^{\times}_J(F) \subseteq \bigcup_{\substack{P\sim \Sigma\\ P \in \euscr{S}_{\D}}} P.$$
\end{theorem}

\begin{remark}
Note that while $S$ is assumed to be upper triangular, the D-special varieties that we intersect with $V_{\bar{c}}$ to get atypical components are arbitrary. In particular, we can choose $S=C^{3n}$ since it is upper triangular D-special. 

\end{remark}

\begin{proof}[Proof of Theorem \ref{DMZPD}]
Fix $\bar{c}$ and consider the variety $V_{\bar{c}}$. Let $X \subseteq V_{\bar{c}} \cap R$ be an atypical component in $S$ where $R\subseteq S$ is D-special. Let also $C_0 \subseteq C$ be a finitely generated subfield of $C$ over which $V_{\bar{c}}, R$ and $X$ are defined. Pick a tuple $\bar{J} \in X(F) \cap E_J^{\times}(F)$ and consider the differential subfield $K:= C_0 \langle \bar{J} \rangle$ of $F$ generated by $C_0$ and $\bar{J}$. By the Seidenberg embedding theorem, $K$ can be embedded into the field of meromorphic functions (of one variable) over some complex domain $W$. Then $\bar{J} = J_{\bar{g}}(\bar{z})$ for some $\bar{g}\in \GL_2(\mathbb{C})$ where $z_i = z_i(w)$ is an analytic function of $w$ defined on $W$. This implies, in particular, that $R$ is $J_{\bar{g}}$-special. Let $Y\subseteq \mathbb{C}^{3n}$ be the complex locus of $J_{\bar{g}}(\bar{z}(w))$. Then it is contained in an analytic component $A$ of the intersection $X \cap J_{\bar{g}}(\mathbb{H}^{\bar{g}})$. In particular, $A$ has no constant coordinates. By Proposition \ref{prop-WMZPD1}, there is a finite collection $\Sigma$ of proper $j$-special subvarieties of $T$, depending only on $V$, such that $\pr_j A$ is contained in some $T' \in \Sigma$. Therefore $\bar{j}\in T'$ and by Lemma \ref{lemma-over-C} $\bar{J}$ is contained in a D-special variety defined over $C$ and associated with $T'$. 
This finishes the proof.
\end{proof}

\section{Functional Modular Zilber--Pink with Derivatives}\label{section-WMZPD}
\setcounter{equation}{0}
\setcounter{claim}{0}

\subsection{D-broad varieties}\label{EC}

In \cite{Aslanyan-Eterovic-Kirby-Diff-EC-j} Aslanyan, Eterovi\'c and Kirby  proved some \emph{Existential Closedness} (henceforth referred to as EC) statements which show roughly that if for a system of equations involving the relation $E_{(z,J)}$ having a solution does not contradict Ax--Schanuel then there is a solution in a differentially closed field. One of those results (in fact, its proof)  will be used in the proof of a functional version of the MZPD conjecture, so we discuss it and some related results here. The reader is referred to \cite{Aslanyan-adequate-predim,Aslanyan-Eterovic-Kirby-Diff-EC-j} for details.

In this subsection all differential fields are ordinary, i.e. have only one derivation.

\begin{notation}
Let $F$ be a field, $n$ be a positive integer, $k \leq n$ and $1\leq i_1 < \ldots < i_k \leq n$. For $\bar{i}=(i_1,\ldots,i_k)$ define the following projection maps. 
\begin{itemize}[leftmargin = 0.5cm]
    \item $\pr_{\bar{i}}:F^{n} \rightarrow F^{k}$ is the map
$$\pr_{\bar{i}}:(x_1,\ldots,x_n)\mapsto (x_{i_1},\ldots,x_{i_k}).$$

\item $\pi_{\bar{i}}:F^{4n}\rightarrow F^{4k}$ is defined by
$$\pi_{\bar{i}}:(\bar{x},\bar{y},\bar{y}',\bar{y}'')\mapsto (\pr_{\bar{i}}\bar{x},\pr_{\bar{i}}\bar{y},\pr_{\bar{i}}\bar{y}',\pr_{\bar{i}}\bar{y}'').$$

\item $\ppr_{\bar{i}}:F^{3n}\rightarrow F^{3k}$ denotes the map
$$\ppr_{\bar{i}}:(\bar{y},\bar{y}',\bar{y}'')\mapsto (\pr_{\bar{i}}\bar{y},\pr_{\bar{i}}\bar{y}',\pr_{\bar{i}}\bar{y}'').$$

\item $\pr_{\bar{y}}:F^{3n}\rar F^n$ and $\ppr_{\bar{y}}:F^{4n}\rar F^n$ are the maps $$\pr_{\bar{y}}:  (\bar{y},\bar{y}',\bar{y}'')\mapsto \bar{y} \mbox{ and } \ppr_{\bar{y}}:  (\bar{z},\bar{y},\bar{y}',\bar{y}'')\mapsto \bar{y},$$ i.e. they are the projections $\pr_j$ and $\ppr_j$ considered in the previous sections.
\end{itemize}
\end{notation}

\begin{definition}
Let $F$ be an algebraically closed field. An irreducible algebraic variety $V \subseteq F^{4n}$ is $J$-\emph{broad} if for any $1\leq i_1 < \ldots < i_k \leq n$ we have $\dim \pi_{\bar{i}} (V) \geq 3k$. We say $V$ is \emph{strongly $J$-broad} if the strict inequality $\dim \pi_{\bar{i}} (V) > 3k$ holds for any $\bar{i}$.
\end{definition}

\begin{definition}
An algebraic variety $V \subseteq F^{4n}$ (or $V \subseteq F^{3n}$) is said to be $J$-\emph{free} if $\ppr_{\bar{y}} V$ (respectively $\pr_{\bar{y}} V$) is not contained in a proper $j$-special subvariety of $F^{n}$.
\end{definition}


\begin{theorem}[EC, {\cite[Theorem 3.6]{Aslanyan-Eterovic-Kirby-Diff-EC-j}}]\label{theorem-EC}
Let $(F; +, \cdot, D)$ be a differentially closed field and let $V\subseteq F^{4n}$ be a strongly $J$-broad and $J$-free variety defined over the field of constants $C$. Then $V(F) \cap E_{(z,J)}^{\times}(F) \neq \emptyset$.
\end{theorem}

We will need this theorem in the next subsection, but it is in fact more convenient to use the method of its proof rather than its statement. So we will not refer to this theorem, but we could actually use its statement instead (which was done in earlier versions of the paper), and the two approaches are equivalent.

\begin{definition}\label{Def-D-broad}
Let $V \subseteq F^{3n}$ be an irreducible variety with a D-special closure $S$ and let $T$ be the $j$-special closure of $\pr_{\bar{y}} V$. Then $V$ is said to be \emph{strongly D-broad} if for all $1\leq i_1 < \ldots < i_k \leq n$ 
\begin{equation*}
    \dim \ppr_{\bar{i}} V > \dim \ppr_{\bar{i}} S - \dim \pr_{\bar{i}} T.
\end{equation*}
\end{definition}

\begin{remarks}
\begin{itemize}[leftmargin = 0.5cm]
    \item[]
    \item It is easy to see that if $V$ is strongly D-broad with respect to a particular choice of $S$ then its $j$-blocks are D-special. So $S$ is the unique D-special closure of $V$ which is simply the product of its $j$-blocks, and $T=\pr_{\bar{y}} S$. In particular, strong D-broadness does not depend on the choice of $S$. 

    \item Theorem \ref{theorem-EC} can be used to show that in a differentially closed field D-broad varieties defined over the field of constants contain $E_J^{\times}$-points. In $\aleph_0$-saturated differentially closed fields they also contain generic $E_J^{\times}$-points. 

    Weak Ax--Schanuel (Theorem \ref{weak-Ax}) implies that if $V$ is defined over a field $C$ and there is a point in $V\cap E^{\times}_J$ generic in $V$ over $C$ then $V$ must be strongly D-broad.
    
    \item A $J$-free variety $V\seq F^{3n}$ is strongly D-broad if and only if $F^n\times V$ is strongly $J$-broad.
\end{itemize}
\end{remarks}

\subsection{A Zilber--Pink type statement}

\begin{definition}
Let $S\subseteq C^{3n}$ be D-special and $V\subseteq S$ be a subvariety. Recall that a \textit{D-atypical} subvariety of $V$ in $S$ is an atypical (in $S$) component $W$ of an intersection $V\cap T$ where $T\subseteq S$ is D-special. If, in addition, $W$ is strongly D-broad then we say that it is \textit{strongly D-atypical}. The \emph{strongly D-atypical} set of $V$ in $S$, denoted $\SAtyp_{\D}(V;S)$, is the union of all strongly D-atypical subvarieties of $V$ in $S$.
\end{definition}

\begin{remark}

One might expect strongly D-atypical subvarieties to be defined as before, that is, if $W$ is D-atypical and no coordinate is constant on $W$ then it is strongly D-atypical. However, the condition of not having any constant coordinates is equivalent to all projections of $W$ having positive dimension. This is actually what the analogue of strong D-broadness would be in the case of $j$ (without derivatives). So, from this point of view, the above notion of strong D-atypicality for D-special varieties is analogous to strong atypicality for $j$-special varieties. Furthermore, broadness and existential closedness have been implicitly used in the proof of weak MZP without derivatives as well, and we did not see those explicitly since the appropriate notion of strong broadness is simpler (equivalent to not having constant coordinates) and the analogue of EC holds trivially there. 

In the following theorem strong D-broadness of atypical subvarieties corresponds to intersection with $E_J^{\times}$ in Theorem \ref{DMZPD}. 


\end{remark}

\begin{theorem}[FMZPD]\label{weak-MZPD}
Let $S\subseteq C^{3n}$ be an upper triangular D-special variety associated with a $j$-special variety $P\subseteq C^n$. Given a parametric family of algebraic subvarieties $(V_{\bar{c}})_{\bar{c} \subseteq C}$ of $S$, there is a finite collection $\Sigma$ of proper $j$-special subvarieties of $P$ such that for every $\bar{c} \subseteq C$, every strongly atypical subvariety of $V_{\bar{c}}$ is contained in a D-special variety associated with some $P'\in \Sigma$. Equivalently,
$$\SAtyp_{\D}(V_{\bar{c}};S) \subseteq \bigcup_{\substack{T\sim \Sigma\\ T \in \euscr{S}_{\D}}} T.$$
\end{theorem}

As pointed out above, Theorem \ref{theorem-EC} can be used to show that in an $\aleph_0$-saturated differentially closed field D-broad varieties defined over the field of constants contain generic $E_J^{\times}$-points (cf. \cite[Theorem 3.8]{Aslanyan-Eterovic-Kirby-Diff-EC-j}). Then the above theorem can be deduced from Theorem \ref{DMZPD} by extending $C$ to an $\aleph_0$-saturated differentially closed field $F$ and working with generic $E_J^{\times}$-points in strongly D-atypical varieties.
Nevertheless, we give a direct differential algebraic proof for Theorem \ref{weak-MZPD}. It has some advantages, for example, it is based on formal properties of the differential equation of the $j$-function and can possibly be adapted to other settings. Note also that Theorem \ref{weak-MZPD} immediately implies Theorem \ref{FMZPD-intro}.

\begin{proof}[Proof of Theorem \ref{weak-MZPD}]
Let $W \subseteq V_{\bar{c}} \cap T$ be a strongly atypical subvariety of $V_{\bar{c}}$ where $T$ is a D-special subvariety of $S$. We know that
$$\dim W > \dim V_{\bar{c}} + \dim T - \dim S.$$

We may assume without loss of generality that $T$ is the D-special closure of $W$ (which is unique since $W$ is strongly D-broad). Indeed, otherwise we could replace $T$ by the D-special closure of $W$ and the above inequality would still hold.\\ 

\noindent \textbf{Step 1.}
Let $\bar{J}:=(\bar{j}, \bar{j}', \bar{j}'') \in W$ be a Zariski generic point over $C$, and let $K:=(C(\bar{J}))^{\alg}$. Let also $d:K\rightarrow \Omega$ be the universal derivation on $K$ over $C$ where $\Omega = \Omega(K/C)$ is the vector space of the abstract differential forms on $K$ over $C$ (see Section \ref{diff-forms}). Then $\dim W = \td(\bar{J}/C) = \dim \Der(K/C) = \dim \Omega(K/C)$.  Consider the differential forms $\omega_i:=\frac{dj_i}{j_i'}-\frac{dj'_i}{j_i''},~ \omega'_i:=\frac{dj'_i}{j_i''}-\frac{dj''_i}{j_i'''}$, where $j_i''' = \eta(j_i,j_i',j_i'')$ is uniquely determined from the equation $\Psi(j_i,j_i',j_i'',j_i''')=0$.  

Let $\Theta:= \Span_K \{ \omega_i, \omega'_i: i=1,\ldots, n \}\subseteq \Omega$, and let $\Lambda:=\Lambda(K/C)$ be the annihilator of $\Theta$, that is,
$$\Lambda=\left\{ D \in \Der(K/C): \frac{Dj_i}{j_i'} = \frac{Dj_i'}{j_i''} = \frac{Dj_i''}{j_i'''},~ i=1,\ldots,n \right\} = \bigcap_i \left( \ker \omega_i \cap \ker \omega_i' \right).$$  Here the annihilator $\Ann(\Theta)$ is a subspace of $(\Omega(K/C))^*$, while $\bigcap_i \ker(\omega_i)$ is a subspace of $\Der(K/C)$. So the above equality makes sense after identifying $(\Omega(K/C))^*$ with its double dual $\Der(K/C)$.
It is clear that $\dim \Lambda = \dim \Omega - \dim \Theta = \dim W - \dim \Theta$.  

Denote the $j$-special variety associated with $T$ by $\tilde{T}$.

\begin{claim}\label{claim-rank}
$\dim \Theta \leq \dim T - \dim \tilde{T}$ and $\dim \Lambda \geq \dim W -(\dim T - \dim \tilde{T}).$
\end{claim}

\begin{proof}

Assume that $j_1$ and $j_2$ are related by a modular equation. Then the coordinates $y_1, y_2$ satisfy the same modular equation on $T$ (since we assumed $T$ is the D-special closure of $W$). Hence $\dim \ppr_{(1,2)} T$ is $3$ or $4$ depending on whether that projection is upper triangular or not. We claim that if it is upper triangular then $\omega_2, \omega'_2\in \Span_K \{ \omega_1, \omega_1' \}$, otherwise $\omega'_2\in \Span_K \{ \omega_1, \omega_1', \omega_2 \}$. 

Let us verify the first assertion first. Suppose $D \in \Der(K/C)$ satisfies $\omega_1(D) = \omega'_1(D) = 0$, that is,
$$ \frac{Dj_1}{j_1'} = \frac{Dj_1'}{j_1''} = \frac{Dj_1''}{j_1'''}.$$ We need to prove that\footnote{If we identify $\Der(K/C)$ with $(\Omega(K/C))^*$ then this means that any linear functional on $\Omega(K/C)$ which vanishes at $\omega_1,\omega_1'$ must also vanish at $\omega_2,\omega_2'$. Hence $\omega_2, \omega'_2\in \Span_K \{ \omega_1, \omega_1' \}$.} $\omega_2(D) = \omega_2'(D) = 0$, i.e.
$$ \frac{Dj_2}{j_2'} = \frac{Dj_2'}{j_2''} = \frac{Dj_2''}{j_2'''}.$$

First observe that if $Dj_1=0$ then $Dj_1'=Dj_1''=Dj_2'=Dj_2''=0$ so there is nothing to prove. Hence we assume $Dj_1\neq 0$. Choose $z_1$ in a differential field extension of $(K;+,\cdot, D)$ such that $Dz_1 = \frac{Dj_1}{j_1'} $. Let $z_2 = g z_1$ and choose $g \in \SL_2(C)$ so that $(z_1,z_2)$ lies in a geodesic variety associated with $\ppr_{(1,2)} T$. Then $$(z_2, j_2, \partial_{z_2} j_2, \partial^2_{z_2} j_2) \in E^{\times}_{(z,J)}$$ and $\partial_{z_2} j_2$ and $j_2'$ satisfy the same algebraic equation over $j_1,j_1',j_1'',j_2$. This implies $j_2' = \partial_{z_2} j_2$ for that equation is linear. Similarly, $j_2'' = \partial^2_{z_2} j_2 = \partial_{z_2} j_2'$ and $j_2''' = \partial^3_{z_2} j_2 = \partial_{z_2} j_2''$. Hence $ \frac{Dj_2}{j_2'} = \frac{Dj_2'}{j_2''} = \frac{Dj_2''}{j_2'''}.$ 

Now assume $\dim \ppr_{(1,2)} T = 4$, that is, $\ppr_{(1,2)} T$ is not upper triangular. Let $U\subseteq C^2$ be a geodesic variety associated with  $\ppr_{(1,2)} T$ and defined by an equation $x_2=\frac{ax_1+b}{cx_1+d}$ with $ad-bc=1$. Assume
\begin{equation*}\label{Dj/j'}
    \frac{Dj_1}{j_1'} = \frac{Dj_1'}{j_1''} = \frac{Dj_1''}{j_1'''}~ \mbox{  and  }~ \frac{Dj_2}{j_2'} = \frac{Dj_2'}{j_2''}.
\end{equation*}

We may assume $Dj_1\neq 0$ as before. We know that $\Phi(j_1,j_2)=0$ for some modular polynomial $\Phi(Y_1,Y_2)$. 
Pick $z_1 \in K$ such that

\begin{equation*}\label{eq-D-special-2}
    \frac{\partial \Phi}{\partial Y_1} (j_1,j_2) \cdot j'_1 + \frac{1}{(cz_1+d)^2}\cdot \frac{\partial \Phi}{\partial Y_2} (j_1,j_2) \cdot j'_2 = 0.
\end{equation*}

The proof of Proposition \ref{prop-S1=S2} shows that
$$Dz_1 = \frac{Dj_1}{j_1'}.$$ Thus, $j_1' = \partial_{z_1}j_1,~ j_1'' = \partial_{z_1}j_1',~ j_1''' = \partial_{z_1}j_1''$ and $(z_1,j_1,j_1',j_1'')\in E_{(z,J)}(K)$. Therefore we can prove as in the upper triangular case that
$$j_2' = \partial_{z_2}j_2,~ j_2'' = \partial_{z_2}j_2',~ j_2''' = \partial_{z_2}j_2''$$ where $z_2=\frac{az_1+b}{cz_1+d}$, which immediately implies the desired equality.

Now if a third coordinate $j_3$ is modularly related to $j_2$ then $j_3'$ is algebraic over $j_1,j_1',j_2,j_2',j_3$ and $j_3''$ is algebraic over $j_1,j_1',j_1'',j_2,j_2',j_2'',j_3,j_3'$ and we can prove as above that $\omega_3, \omega_3' \in \Span_K \{ \omega_1, \omega_1', \omega_2 \}$. In the upper triangular case we obviously would have $\omega_3, \omega_3' \in \Span_K \{ \omega_1, \omega_1' \}$.

Thus, each $j$-block of $T$ of dimension $3$ contributes at most $2$ to $\dim \Theta$ while each $j$-block of dimension $4$ contributes at most $3$. Hence $\dim \Theta \leq \dim T - \dim \tilde{T}$.
\end{proof}

\begin{claim}\label{Lambda-closed-Lie}
$\Lambda$ is closed under the Lie bracket.
\end{claim}
\begin{proof}
Pick two derivations $D_1, D_2 \in \Lambda$ and let $D:=[D_1, D_2]$. Using the equalities
$$\frac{D_ij}{j'} = \frac{D_ij'}{j''} = \frac{D_ij''}{j'''},~ i=1,2,$$ we get
$$\frac{Dj}{j'} = \frac{Dj'}{j''} + \frac{1}{j'} \left( D_2j' \cdot D_1 \frac{j'}{j''} - D_1j' \cdot D_2 \frac{j'}{j''} \right).$$
We claim that the expression in brackets is equal to zero. Indeed, after simplifying it we see that it suffices to prove that 
$$D_2 j' \cdot D_1 j'' = D_1 j' \cdot D_2 j''.$$
This is equivalent to $$\frac{D_1 j'}{D_1 j''} = \frac{D_2 j'}{D_2 j''}.$$ But we know that for $i=1,2$
$$\frac{D_i j'}{D_i j''} = \frac{j''}{j'''}$$ which does not depend on $i$. This shows that $\frac{Dj}{j'} = \frac{Dj'}{j''}.$

Similarly, the equality $\frac{Dj'}{j''} = \frac{Dj''}{j'''}$ follows from $\frac{D_1 j''}{D_1 j'''} = \frac{D_2 j''}{D_2 j'''}.$ We know that $j''' = \eta(j,j',j'')$ where $\eta(Y_1, Y_2, Y_3)$ is a rational function over $\mathbb{Q}$. Hence for $i=1,2$
\begin{align*}
    \frac{D_i j'''}{D_i j''} & = \frac{\partial \eta}{\partial Y_1} (j, j', j'') \cdot \frac{D_ij}{D_ij''} + \frac{\partial \eta}{\partial Y_2} (j, j', j'') \cdot \frac{D_ij'}{D_ij''} + \frac{\partial \eta}{\partial Y_3} (j, j', j'') \cdot \frac{D_ij''}{D_ij''} \\ & = \frac{\partial \eta}{\partial Y_1} (j, j', j'') \cdot \frac{j'}{j'''} + \frac{\partial \eta}{\partial Y_2} (j, j', j'') \cdot \frac{j''}{j'''} + \frac{\partial \eta}{\partial Y_3} (j, j', j'')
\end{align*}
which does not depend on $i$. This finishes the proof. 
\end{proof}

\noindent \textbf{Step 2.}
Let $l:= \dim \Lambda$. By Lemma \ref{Lie-Lemma} we can take a commuting basis $D_1,\ldots, D_l$ of $\Lambda$ and consider the differential field $(K;+,\cdot,D_1,\ldots,D_l)$. Let $\hat{C}:= \bigcap_i\ker(D_i)$ be its field of constants. 



\begin{claim}\label{generic-derivation}
$j_k \notin \hat{C}$ for each $k$ (and hence $j_k', j_k'' \notin \hat{C}$).
\end{claim}
\begin{proof}
This is an adaptation of the proof of \cite[Theorem 3.6]{Aslanyan-Eterovic-Kirby-Diff-EC-j}. 

Assume $\dim W = \dim T - \dim \tilde{T} + r$ where $r\geq 1$ (recall that $W$ is strongly D-broad). Then by Claim \ref{claim-rank} we have $l\geq r$.

Let $J_i:=(j_i,j_i',j_i'')$. It is evident that if one of the coordinates of $J_i$ is in $\hat{C}$ then so are the others (see \cite[Lemma 3.4]{Aslanyan-Eterovic-Kirby-Diff-EC-j}). Hence we may assume that for some $t\geq 1$ none of the coordinates of $J_1, \ldots, J_t$ is in $\hat{C}$ and all coordinates of $J_{t+1},\ldots, J_n$ are in $\hat{C}$. Here $t\neq 0$ for otherwise all derivations $D_i$ would vanish on $K$ which contradicts $l>0$. Assume $t<n$. Let
$$\bar{e}:=(1,\ldots,t),~ \bar{s}:=(t+1,\ldots,n)$$ and $$\bar{v}:= (J_1, \ldots, J_t),~ \bar{w} := (J_{t+1},\ldots,J_n).$$

Since $\hat{C}$ is algebraically closed, there cannot be a modular relation between a coordinate of $\bar{v}$ and a coordinate of $\bar{w}$. Therefore $T = \ppr_{\bar{e}}T \times \ppr_{\bar{s}}T$ and $\tilde{T} = \pr_{\bar{e}}\tilde{T} \times \pr_{\bar{s}}\tilde{T}$. In particular, $\ppr_{\bar{e}} T$ is the D-special closure of $\bar{v}$.

Since $\bar{v}\in E_J^{\times}(K)$, by Theorem \ref{weak-Ax} (weak Ax--Schanuel) we conclude that $$\td(\hat{C}(\bar{v})/\hat{C}) \geq \dim \ppr_{\bar{e}} T - \dim \pr_{\bar{e}}\tilde{T}+\rk(D_i j_k)_{1\leq i\leq l,1\leq k\leq t} = \dim \ppr_{\bar{e}}T - \dim \pr_{\bar{e}}\tilde{T}+ l.$$ 
Further, $$\td(\hat{C}/C)\geq \td(C(\bar{w})/C) = \dim \ppr_{\bar{s}} (W) > \dim \ppr_{\bar{s}}T - \dim \pr_{\bar{s}}\tilde{T} $$ for $W$ is strongly D-broad.

Combining the above inequalities and equalities we get $$\dim W = \td(K/C) = \td(K/\hat{C}) + \td(\hat{C}/C) > \dim T - \dim \tilde{T} + l \geq \dim T - \dim \tilde{T} + r,$$ which is a contradiction. Thus, $t=n$ and no coordinate of $\bar{J}$ is in $\hat{C}$.
\end{proof}

\noindent \textbf{Step 3.} Now let $U_{\bar{d}}\subseteq K^n$ be a $C$-geodesic variety associated with $T$ chosen from the parametric family of all $\hat{C}$-geodesic subvarieties of $K^n$. 
Then $\dim U_{\bar{d}} = \dim \tilde{T}$. 
Extend (if necessary) the differential field $(K; +, \cdot, D_1,\ldots,D_l)$ by adjoining elements $(z_1,\ldots,z_n)\in U_{\bar{d}}$ with 
$$D_s z_i = \frac{D_s j_i}{j'_i},~ i=1,\ldots, n,~ s= 1, \ldots, l.$$ 
Denote the field $K(\bar{z})^{\alg}$ by $F$. Note that if $z_i$ corresponds to a $j$-block of $T$ of dimension $4$ then $z_i$ is algebraic over (and hence belongs to) $K$. In this case we just take $z_i$ to be an element of $K$ satisfying \eqref{eq-D-special-2-0}, and then we will automatically have $D_s z_i = \frac{D_s j_i}{j'_i}$ for all $s$. On the other hand, if $z_i$ corresponds to a $j$-block of $T$ of dimension $3$ then it is transcendental over $K$. In this case we choose $z_i$'s to be algebraically independent over $K$ if they correspond to different $j$-blocks of dimension $3$ (and they must be linked by $\SL_2(C)$-relations if they correspond to the same $j$-block). Then we extend $D_s$ by defining $D_s z_i = \frac{D_s j_i}{j'_i}$. This makes $F$ a differential field extension of $K$ with derivations $D_1,\ldots,D_l$, and clearly $(\bar{z}, \bar{j},\bar{j}',\bar{j}'')\in E^{\times}_{(z,J)}(F)$. Simple calculations as in Claim \ref{Lambda-closed-Lie} show that the derivations commute on $z_i$'s and hence on $F$. 

The field of constants of $F$ may be larger than $\hat{C}$ but by abuse of notation we still denote it by $\hat{C}$. Let us also stress that here we choose $\bar{z}$ from $U_{\bar{d}}$, which means that the geodesic relations linking $z_i$'s come from $\SL_2(C)$ rather than $\SL_2(\hat{C})$.   Further, as we saw in Section \ref{diff-forms}, $\rk \Jac (\bar{z}) = \rk \Jac (\bar{J}) = \dim \Lambda$. Thus
\begin{gather*}
    \dim U_{\bar{d}} \times V_{\bar{c}} = \dim U_{\bar{d}} + \dim V_{\bar{c}} < \dim \tilde{T}+ \dim W +\dim S - \dim T =\\
    = 3\dim P + \dim W - (\dim T - \dim \tilde{T}) \leq 3\dim P + \rk \Jac (\bar{z}).
\end{gather*}

Now we apply the uniform Ax--Schanuel with derivatives to the parametric family $(U_{\bar{d}} \times V_{\bar{c}})_{\bar{c},\bar{d}\subseteq \hat{C}}$ and get a finite collection $\Sigma$ of proper $j$-special varieties of $P$, depending on this parametric family only (which, in turn, depends only on $V$ and is independent of $T$ and $W$), such that $\bar{j}\in P'$ for some $P' \in \Sigma$. Then $\bar{J}$ belongs to a D-special variety associated with $P'$ which is not necessarily defined over $C$ (since it is possible that $C\subsetneq \hat{C}$). However, since $W$ and $T$ are defined over $C$ and $W$ is strongly D-broad, it follows that there is a D-special variety $S'$ associated with $P'$ and defined over $C$ such that $\bar{J}\in S'$. So we conclude that $W\subseteq S'$ as $\bar{J}$ is generic in $W$ over $C$. 
\end{proof}

\subsection{A differential algebraic proof of DMZPD}

Now we prove that FMZPD implies DMZPD, which gives a differential algebraic proof for DMZPD (modulo the proof of Ax-Schanuel, which is not differential algebraic). 

\begin{lemma}\label{lemma-atypical-subvariety}
Assume $X, Y \subseteq Z$ are algebraic varieties, with $X$ irreducible, and $\bar{x}$ is a smooth point of $X$ which belongs to an atypical component of the intersection $X\cap Y$ in $Z$. Then for any irreducible subvariety $X'\subseteq X$ containing $\bar{x}$ the intersection $X' \cap Y$ is atypical in $Z$ and $\bar{x}$ belongs to an atypical component of that intersection.
\end{lemma}
\begin{proof}

Let $A\subseteq X\cap Y$ be an atypical component with $\bar{x}\in A$. Choose a component $A'$ of $A\cap X'$ containing $\bar{x}$. Then by Theorem \ref{dimension-intersect-nonsmooth} (for algebraic varieties) we get
\begin{gather*}
    \dim A' \geq \dim A+\dim X' - \dim X >\\
     (\dim X + \dim Y - \dim Z) + \dim X' - \dim X = \dim X' + \dim Y - \dim Z. \qedhere
\end{gather*}
\end{proof}

\begin{proposition}\label{JWEC-WMZPD-j'}
\textup{FMZPD} implies \textup{DMZPD}. 
\end{proposition}

\begin{proof}
Let $T'$ be a D-special variety which intersects $V_{\bar{c}}$ atypically. Assume $\bar{J}=(\bar{j},\bar{j}',\bar{j}'')\in V_{\bar{c}}(F) \cap T'(F) \cap E^{\times}_J(F)$ belongs to an atypical component of $V_{\bar{c}}\cap T'$. The point $\bar{J}$ has a unique D-special closure which we denote by $T$. We claim that $\bar{J}$ belongs to an atypical component of $V_{\bar{c}}\cap T$ in $S$.

Although D-special varieties are not smooth in general (and even  modular curves have singularities), $E_J^{\times}$-points on D-special varieties are non-singular since such points are generic in $j$-blocks (Corollary \ref{cor-D-special-EJ-generic}) and a D-special variety is just the product of its $j$-blocks. Thus, $\bar{J}$ is a smooth point of $T'$, hence by Lemma \ref{lemma-atypical-subvariety} $\bar{J}$ belongs to an atypical component $W$ of the intersection $V_{\bar{c}}\cap T$ in $S$. Then $T$ is the unique D-special closure of $W$, and $\tilde{T}:= \pr_{\bar{y}} T$ is the $j$-special closure of $\pr_{\bar{y}} W$.

Now Theorem \ref{weak-Ax} implies $$\dim W \geq \td_CC(\bar{J}) > \dim T - \dim \tilde{T}.$$ Moreover, this inequality holds for all projections of $W$ and hence it is strongly D-broad. Therefore the hypotheses of Theorem \ref{weak-MZPD} are satisfied, and we are done.
\end{proof}

\section{Functional Modular Andr\'{e}--Oort with Derivatives}\label{section-FMAOD}
\setcounter{equation}{0}

The following Modular Andr\'{e}--Oort with Derivatives conjecture is a special case of MZPD and was proposed by Pila in his aforementioned unpublished notes (see also \cite{Spence}).

\begin{conjecture}[MAOD]
For every algebraic variety $V \subsetneq \mathbb{C}^{3n}$ there is a finite collection $\Sigma$ of proper $\mathbb{H}$-special subvarieties of $\mathbb{H}^n$ such that every $J$-special subvariety of $V$ is contained in a $J$-special variety of the form $\langle \langle \bar{\gamma} U \rangle \rangle$ for some $\bar{\gamma}\in \SL_2(\mathbb{Z})^n$ and some $U \in \Sigma$.
\end{conjecture}

Note that here one does not intersect $J$-special subvarieties with the image of $J$ since these varieties always contain $J$-points. We prove a functional analogue of this conjecture.

\begin{theorem}[FMAOD]\label{weak-MAOD}
Let $S\subseteq C^{3n}$ be an upper triangular D-special variety associated with a $j$-special variety $T\subseteq C^n$. Given a parametric family of algebraic subvarieties $(V_{\bar{c}})_{\bar{c} \subseteq C}$ of $S$, there is a finite collection $\Sigma$ of proper $j$-special subvarieties of $T$ such that for every $\bar{c} \subseteq C$, if $V_{\bar{c}}\subsetneq S$ then every D-special subvariety of $V_{\bar{c}}$ is contained in a D-special variety associated with some $T'\in \Sigma$.
\end{theorem}

\begin{lemma}
D-special varieties are strongly D-broad. In particular, if $V\subsetneq S\subseteq C^{3n}$ is a proper subvariety of a D-special variety $S$ and $T\subseteq V$ is D-special then $T$ is a strongly D-atypical subvariety of $V$.
\end{lemma}
\begin{proof}
Straightforward.
\end{proof}

\begin{proof}[Proof of Theorem \ref{weak-MAOD}]
This follows from Theorem \ref{weak-MZPD} along with the above lemma.
\end{proof}

\begin{example}[cf. \cite{Spence}, example after Definition 1.5]\label{example-AO-infinite}
Let $V\subseteq C^6$ and $T\subseteq C^2$ be defined by an equation $\Phi(y_1,y_2)=0$ where $\Phi$ is a modular polynomial, i.e. $T = \pr_{\bar{y}} V$. Then all D-special varieties associated with $T$ are contained in $V$ and they are maximal in $V$. Thus, in Theorem \ref{weak-MAOD} even a single variety $V$ may contain a whole parametric family of D-special varieties (associated with finitely many $j$-special varieties).

Here $V$ is not $J$-free. To get a similar example with a $J$-free variety let $W\seq C^9$ be defined by a single equation $\Phi_N(y_1,y_2)+\Phi_M(y_2,y_3) = 0$ for some $N, M$. Now if $T\seq C^3$ is a component of the variety defined by $\Phi_N(y_1,y_2) = \Phi_M(y_2,y_3) = 0$ then all D-special varieties associated with $T$ are contained in $V$, and they are maximal in $V$.
\end{example}


\addtocontents{toc}{\protect\setcounter{tocdepth}{1}}
\subsection*{Acknowledgements} I am grateful to Jonathan Pila and Haden Spence for useful discussions on the Modular Zilber--Pink with Derivatives and Andr\'{e}--Oort with Derivatives conjectures. I would also like to thank the referee for reading the paper and for the positive feedback.

\addtocontents{toc}{\protect\setcounter{tocdepth}{2}}

\bibliographystyle {alpha}
\bibliography {ref}

\end{document}